\documentclass[12pt]{amsart}
\usepackage{latexsym,amsmath,amsthm,amssymb,amscd}
\newtheorem{thm}{Theorem}[section]
\newtheorem{lem}[thm]{Lemma}
\newtheorem{corollary}[thm]{Corollary}
\newtheorem{remark}[thm]{Remark}

\begin{document}
\renewcommand{\d}{\textrm{d}}
\def\Ad{{\rm Ad}}
\def\rad{{\rm rad}}
\def\val{{\rm val}}
\def\diag{{\rm diag}}
\def\End{{\rm End}}
\def\Fr{{\rm Fr}}
\def\Gal{{\rm Gal}}
\def\GL{{\rm GL}}
\def\I{{\rm I}}
\def\Id{{\rm Id}}
\def\rs{{\rm rs}}
\def\fa{{\mathfrak a}}
\def\fc{{\mathfrak c}}
\def\haut{{\rm ht}}
\def\ldee{{\ell(\gamma)}}
\def\imag{{\rm im}}
\def\Im{{\rm Im}}
\def\Mat{{\rm M}}
\def\Nrd{{\rm Nrd}}
\def\integers{{\mathfrak o}}
\def\fR{{\mathfrak R}}
\def\Lie{{\rm{Lie}}}
\def\LieA{{\mathfrak a}}
\def\LieN{{\mathfrak n}}
\def\re{{\rm Re}}
\def\Rep{{\rm Rep}}
\def\Res{{\rm Res}}
\def\WDRep{{\rm WDRep}}
\def\rec{{\rm rec}}
\def\Hom{{\rm Hom}}
\def\Ind{{\rm Ind}}
\def\cInd{{\rm c\!\!-\!\!Ind}}
\def\pr{{\rm pr}}
\def\supp{{\rm supp}}
\def\Mod{{\rm Mod}}
\def\ii{{\rm i}}
\def\rr{{\rm r}}
\def\vol{{\mu}}
\def\Pbar{{\rm \overline P}}
\def\Nbar{{\rm \overline N}}
\def\St{{\rm St}}
\def\Stab{{\rm Stab}}
\def\assoc{{assoc}}
\def\sw{{\rm sw}}
\def\Sw{{\rm Sw}}
\def\Ar{{\rm A}}
\def\Sym{{\rm Sym}}
\def\Cusp{{\rm Cusp}}
\def\Irr{{\rm Irr}}
\def\Irrt{{{\rm Irr}^{\rm t}}}
\def\Psit{{{\Psi}^{\rm t}}}
\def\temp{{{\rm t}}}
\def\trace{{{\rm trace}}}
\def\simple{{{\rm s}}}
\def\Irru{{{\rm Irr}^{\rm u}}}
\def\mes{{\rm mes}}
\def\Tor{{\rm Tor}}
\def\Nor{{\rm N}}
\def\ie{{\it i.e.,\,}}
\def\tr{{\rm tr}}
\def\Omegat{{{\Omega}^{\rm t}}}
\def\cD{{\mathcal D}}
\def\cE{{\rm E}}
\def\cEess{{\rm E}_{{\rm ess}}}
\def\JL{{\rm{JL}}}
\def\cH{{\mathcal H}}
\def\cL{{\mathcal L}}
\def\cK{{\mathcal K}}
\def\spin{{\rm sp}}
\def\SL{{\rm SL}}
\def\Sp{{\rm Sp}}
\def\SU{{\rm SU}}
\def\red{{\rm nd}}
\def\CC{{\mathbb C}}
\def\QQ{{\mathbb Q}}
\def\RR{{\mathbb R}}
\def\ZZ{{\mathbb Z}}
\def\Afr{{\mathfrak A}}
\def\Bfr{{\mathfrak B}}
\def\Cfr{{\mathfrak C}}
\def\Ofr{{\mathfrak O}}
\def\Pfr{{\mathfrak P}}

\title{Plancherel measure for $\GL(n,F)$ and $\GL(m,D)$: explicit formulas and Bernstein decomposition}
\author{Anne-Marie Aubert and Roger Plymen}

\date{}

\maketitle

\pagestyle{myheadings}
\markboth{A.-M.~Aubert and R.~Plymen}{Plancherel measure for $\GL(n)$:
explicit formulas and Bernstein decomposition}
\begin{abstract}  Let $F$ be a nonarchimedean local field, let $D$ be a division algebra over $F$,
let $\GL(n) = \GL(n,F)$.  Let $\nu$ denote Plancherel measure for
$\GL(n)$.  Let $\Omega$ be a component in the Bernstein variety
$\Omega(\GL(n))$. Then $\Omega$ yields its fundamental invariants:
the cardinality $q$ of the residue field of $F$, the sizes $m_1,
\ldots, m_t$, exponents $e_1, \ldots, e_t$, torsion numbers $r_1,
\ldots, r_t$, formal degrees $d_1, \ldots, d_t$ and conductors
$f_{11}, \ldots, f_{tt}$.  We provide explicit formulas for the
Bernstein component $\nu_{\Omega}$ of Plancherel measure in terms
of the fundamental invariants.   We prove a transfer-of-measure
formula for $\GL(n)$ and establish some new formal degree
formulas. We derive, via the Jacquet-Langlands correspondence, the
explicit Plancherel formula for $\GL(m,D)$.


\end{abstract}
\maketitle
\section{Introduction}

In this article we provide an explicit Plancherel formula for the
$p$-adic group $\GL(n)$.   Moreover, we determine explicitly the
Bernstein decomposition of Plancherel measure, including all
numerical constants.

Let $F$ be a nonarchimedean local field with ring of integers
$\integers_F$, let $G = \GL(n) =
\GL(n,F)$.  We will use the standard normalization of Haar measure
on $\GL(n)$ for which the volume of $\GL(n, \integers_F)$ is $1$.
Plancherel measure $\nu$ is then uniquely determined by the
equation
\[
f(g) = \int \trace\,\pi(\lambda(g)f^{\vee})d\nu(\pi)
\]
for all $g \in G, f \in \mathcal{C}(G)$, where $f^{\vee}(g) =
f(g^{-1})$.

The Harish-Chandra Plancherel Theorem expresses the Plancherel
measure in the following form:
\[d\nu(\omega)\, =\,
c(G|M)^{-2}\,\gamma(G|M)^{-1}\, \mu_{G|M}(\omega)\, d(\omega)\, d\omega \]
where $M$ is a Levi subgroup of $G$, $\omega \in \cE_2(M)$ the
discrete series of $M$, $c(G|M)$ and $\gamma(G|M)$ are certain
constants, $\mu_{G|M}$ is a certain rational function, $d(\omega)$ is
the formal degree of $\omega$, and $d\omega$ is the Harish-Chandra
canonical measure.

In this article we determine explicitly
\[
c(G|M)^{-2}\gamma(G|M)^{-1} \mu_{G|M}(\omega)\, d(\omega)\,
d\omega
\] for $\GL(n)$.

The support of Plancherel measure $\nu$ admits a Bernstein
decomposition \cite{P2} and therefore $\nu$ admits a canonical
decomposition
\[ \nu = \bigsqcup \nu_{\Omega} \] where $\Omega$ is a component
in the Bernstein variety $\Omega(G)$.  We determine explicitly the
Bernstein component $\nu_{\Omega}$ for $\GL(n)$.

We can think of $\Omega$ as a vector of irreducible supercuspidal
representations of smaller general linear groups.
If the vector is \[(\sigma_1, \ldots, \sigma_1, \ldots, \sigma_t,
\ldots, \sigma_t)\] with $\sigma_i$ repeated $e_i$ times, $1 \leq
i \leq t$, and $\sigma_1$, $\ldots$, $\sigma_t$ pairwise distinct
(after unramified twist) then we say that $\Omega$ has {\it
exponents} $e_1$, $\ldots$, $e_t$.

Each representation $\sigma_i$ of $\GL(m_i)$ has a {\it torsion
number}: the order of the cyclic group of all those unramified
characters $\eta$ for which $\sigma_i \otimes \eta \cong
\sigma_i$.  The torsion number of $\sigma_i$ will be denoted $r_i
$.

We may choose each representation $\sigma_i$ of $\GL(m_i)$ to be
unitary: in that case $\sigma_i$ has a {\it formal degree} $d_i =
d(\sigma_i)$.  We have $0 < d_i < \infty$.

We will denote by $f_{ij} = f(\sigma_i^{\vee}\times\sigma_j)$ the
conductor of the pair $\sigma_i^{\vee}\times\sigma_j$.  An
explicit conductor formula is obtained in the article by Bushnell,
Henniart and Kutzko \cite{BHKcond}.

In this way, the Bernstein component $\Omega$ yields up the
following {\it fundamental invariants}:
\begin{itemize}
\item the cardinality $q$ of the residue field of $F$ \item the
sizes $m_1, m_2, \ldots, m_t$ of the smaller general linear groups
\item the exponents $e_1, e_2, \ldots, e_t$ \item the torsion
numbers $r_1, r_2, \dots, r_t$ \item the formal degrees $d_1, d_2,
\ldots, d_t$ \item the conductors for pairs $f_{11}, f_{12},
\ldots, f_{tt}$.
\end{itemize}

Our Plancherel formulas are built from precisely these numerical
invariants.

If $\Omega$ has the single exponent $e$, then the fundamental
invariants yielded up by $\Omega$ are $q, m, e, r, d, f$.   The
component $\Omega$ determines a representation in the discrete
series of $\GL(n)$, namely the generalized Steinberg
representation $\St(\sigma,e)$. The formal degree of $\pi =
\St(\sigma,e)$ is given by the following new formula, which is
intricate, but depends only on the fundamental invariants of
$\Omega$, in line with our general philosophy:
\[\frac{d(\pi)}{d(\sigma)^e}=\frac{m^{e-1}}{r^{e-1}e}\,\cdot\,
q^{(e^2-e)(f(\sigma^\vee\times\sigma)+r-2m^2)/2} \,\cdot\,
\frac{(q^r-1)^e}{q^{er}-1}
\,\cdot\,\frac{|\GL(em,q)|}{|\GL(m,q)|^e}.\]

In section 2, we give a pr\'{e}cis of the background material
which we need, following the recent article of Waldspurger
\cite{Wald}.

The Langlands-Shahidi formula gives the rational function
$\mu_{G|M}$ as a ratio of certain $L$-factors and
$\epsilon$-factors \cite{Sh1}. In sections 3--4 we compute
explicitly the expression
\[
c(G|M)^{-2}\gamma(G|M)^{-1} \mu_{G|M}(\omega) d\omega \] when $M$
is a maximal parabolic.  The resulting formula is stated in
Theorem 4.4: in this formula we correct certain misprints in
\cite[p. 292 -- 293]{Sh2}.

In section 5, we compute the Plancherel density $\mu_{G|M}$ in the
general case by using the Harish-Chandra product formula and we
give the explicit Bernstein decomposition of Plancherel measure.

As a special case, we derive the explicit Plancherel formula for
the (extended) affine Hecke algebra $\mathcal{H}(n,q)$.

We have, in effect, extended the classical formula of Macdonald
\cite{M1}, \cite[Theorem 5.1.2]{M2} from the spherical component
of $\GL(n)$ to the whole of the tempered dual.


The Plancherel formulas for $\GL(n,F)$ and $\GL(m,D)$ are
dominated by \emph{repeating patterns}, which we now attempt to
explain.  The repeating patterns are expressed by
transfer-of-measure theorems, of which the first is as follows.
With $j = 1,2$, let $F_j$ be a nonarchimedean local field and let
$\Omega_j$ be a component in the Bernstein variety of $\GL(n_j,
F_j)$. Let $\nu^{(j)}$ denote the Plancherel measure of $\GL(n_j,
F_j)$.  If $\Omega_1, \Omega_2$ share the same fundamental
invariants, then
\[
\nu^{(1)}_{\Omega_1} = \nu^{(2)}_{\Omega_2}.\]

The next transfer-of-measure theorem is more surprising. Let
$\Omega$ be a component in the Bernstein variety of $\GL(n,F)$,
and let $\nu$ be Plancherel measure.   Let $\Omega$ have the
fundamental invariants $(q,m,e,r,d,f)$. Let $K/F$ be an extension
field with $q_K = q^r$.  Let $G_0: = \GL(e,K)$, let $\Omega_0$ be
a component in the Bernstein variety of $G_0$, and let $\nu^{(0)}$
be Plancherel measure. If $\Omega_0$ has fundamental invariants
$(q^r,1,e,1,1,1)$ then $\nu_{\Omega}$ and $\nu^{(0)}_{\Omega_0}$
are \emph{proportional}, \ie
\[\nu_{\Omega} = \kappa \cdot \nu^{(0)}_{\Omega_0}\]
where $\kappa = \kappa(q,m,e,r,d,f)$.  This phenomenon was first
noted by Bushnell, Henniart, Kutzko \cite[Theorem 4.1]{BHK}, working
in the context of types and Hilbert algebras.  We reconcile our
result for $\GL(n)$ with (a special case of) their result by
proving that
\[
\kappa(q,m,e,r,d,f) = vol(J)^{-1} \cdot vol(I_0) \cdot
\dim(\lambda)\] where $(J,\lambda)$ is an $\Omega$-type, $I_0$ is
an Iwahori subgroup of $G_0$: for this result, see Theorem 6.12.
Theorem 5.7, which in essence is the Harish-Chandra product
formula, then allows one to compute the Plancherel measure
$\nu_{\Omega}$ for any component $\Omega$.


Using the explicit value for the formal degree of any
representation in the discrete series of $G$ previously obtained
by Silberger and Zink, we show that the comparison formula between
formal degrees, proved by Corwin, Moy, Sally in the tame case
\cite{CMS}, is valid in general.

In the last section of the paper we consider the case of a group
$\GL(n',D)$, where $D$ is a central division algebra of index $d$
over over $F$. We extend the transfer-of-measure result of Arthur
and Clozel \cite[pp.~88-90]{AC} to the case when $F$ is of
positive characteristic, by using results of Badulescu.

Let $G' = \GL(n',D)$, $G = \GL(n,F)$ with $n = dn'$.
Let $\nu', \nu$ denote the Plancherel measure for $G',G$, each with
the standard normalization of Haar measure on $G',G$.
Let $\JL\colon\cE_2(G')\to\cE_2(G)$ denote the Jacquet-Langlands correspondence.
Then we have
\[
d\nu'(\omega') = \lambda(D/F) \cdot d\nu(\JL(\omega'))
\]
where
\[
\lambda(D/F) = \prod(q^m - 1)^{-1}\] the product taken over all
$m$ such that $1 \leq m \leq n-1, m \neq 0 \, \mod \, d$.

For example, let $G' = \GL(3,D), G = \GL(6,F)$ with $D$ of index
$2$.   Then we have
\[
d\nu'(\omega') = (q-1)^{-1}(q^3-1)^{-1}(q^5-1)^{-1} \cdot
d\nu(\JL(\omega'))
\]

Our proof of this is in local harmonic analysis, cf \cite[p. 88 -
90]{AC}.

\emph{Historical Note.} The Harish-Chandra Plancherel Theorem, and
the Product Theorem for Plancherel Measure, were published
posthumously in his collected papers in 1984, see \cite{HC}. The
theorems were stated without proof (although Harish-Chandra had
apparently written out the proofs). At this point, we quote from
Silberger's article \cite{Si}, published in 1996:
\begin{quotation} In \cite{HC} Harish-Chandra has summarized the
theory underlying the Plancherel formula for G and sketched a
proof of the Plancherel theorem. To complete this sketch it seems
to this writer that details need to be supplied justifying only
one assertion of \cite{HC}, namely Theorem 11. Every other
assertion in this paper can be readily proved either by using
prior published work of Harish-Chandra or the present author's
notes on Harish-Chandra's lectures.\end{quotation}  For
Silberger's Notes, published in 1979, see \cite{Si1}.  Complete
and detailed proofs were finally published by Waldspurger in 2003,
see \cite[V.2.1, VIII.1.1]{Wald}. None of these sources contains
any explicit computations for $\GL(n)$.

We thank the referee for several valuable comments, and especially
for supplying an idea which led to a short proof of Theorem 3.1.
We thank Ioan Badulescu, Volker Heiermann, and Wilhelm Zink for
valuable conversations on the central simple algebras case. Some
of the results in this article have been announced in \cite{AP}.

\section{The Plancherel Formula after Harish-Chandra}

We shall follow very closely the notation and terminology in
Waldspurger \cite{Wald}.

Let $\cK =\GL(n,\integers_F)$. Let $H$ be a closed subgroup of
$G=\GL(n,F)$. We use the {\it standard} normalization of Haar
measures, following \cite[I.1, p.240]{Wald}. Then Haar measure
$\vol_H$ on $H$ is chosen so that $\vol_H(H \cap \cK) = 1$. If $Z
= A_G$ is the centre of $G$ then we have $\vol_Z(Z \cap \cK) = 1$.
If $H = G$ then Haar measure $\vol = \vol_G$ is normalized so that
the volume of $\cK$ is $1$.

Denote by $\Theta$ the set of pairs $(\mathcal{O}, P = MU)$ where
$P$ is a semi-standard parabolic subgroup of $G$ and $\mathcal{O}
\subset \cE_2(M)$ is an orbit under the action of $\Im\,X(M)$.
(Here $\cE_2(M)$ is the set of equivalence classes of the discrete
series of the Levi subgroup $M$, and $\Im\,X(M)$ is the group of
the unitary unramified characters of $M$.)

Two elements $(\mathcal{O}, P = MU)$ and $(\mathcal{O'}, P' =
M'U')$ are {\it associated} if there exists $s \in W^G$ such that
$s \cdot M = M', s\mathcal{O} = \mathcal{O}'$.   We fix a set
$\Theta/\assoc$ of representatives in $\Theta$ for the classes of
association.  For $(\mathcal{O}, P = MU) \in \Theta$, we set
$W(G|M)=\{s\in W^G: s\cdot M=M\}/W^M$, and
\[\Stab(\mathcal{O}, M) = \{s \in W(G|M): s \mathcal{O} =
\mathcal{O}\}.
\]

Let $\mathcal{C}(G)$ denote the Harish-Chandra Schwartz space of
$G$ and let $I_P^G\omega$ denote the normalized induced representation from
$\omega$.  Let $f \in \mathcal{C}(G)$, $\omega \in
\cE_2(M)$. We will write
\[
\pi = I^G_P \omega,\;\;\; \pi(f) = \int f(g)\pi(g)dg,\;\;\;
\theta^G_{\omega}(f) = \trace\, \pi(f).
\]

\begin{thm} \label{HCformula}
{\rm The Plancherel Formula \cite[VIII.1.1]{Wald}.} For each $f \in
\mathcal{C}(G)$ and each $g \in G$ we have
\[
f(g) = \sum c(G|M)^{-2}\gamma(G|M)^{-1}|\Stab(\mathcal{O}, M)|^{-1}
\int_{\mathcal{O}}\mu_{G|M}(\omega)d(\omega)\theta^G_{\omega}(\lambda(g)f^{\vee})d\omega
\]
where the sum is over all the pairs $(\mathcal{O}, P = MU) \in
\Theta/\assoc$.
\end{thm}

\smallskip

Note that
\begin{equation} \label{jfunction}
\mu_{G|M}(\omega) \cdot c(G|M)^{-2} \cdot
\gamma(G|M)^{-1}=\gamma(G|M)\,\cdot\,j(\omega)^{-1},
\end{equation}
where $j$ denotes the composition of intertwining operators defined
in \cite[IV.3~(2)]{Wald}.

\smallskip

The map \[(\mathcal{O},P = MU) \to \Irrt(G),\; \omega \mapsto
I_P^G \omega\] determines a {\it bijection}
\[
\bigsqcup (\mathcal{O}, P = MU)/\Stab(\mathcal{O},M)
\longrightarrow \Irrt(G).
\]

The tempered dual $\Irrt(G)$ acquires, by transport of structure,
the structure of {\it disjoint union of countably many compact
orbifolds}.

According to \cite[V.2.1]{Wald}, the function $\mu_{G|M}$ is a
rational function on $\mathcal{O}$.  We have $\mu_{G|M}(\omega)
\geq 0$ and $\mu_{G|M}(s\omega) = \mu(\omega)$ for each $s \in
W^G, \omega \in \mathcal{O}$.   This invariance property implies
that $\mu$ {\it descends} to a function on the orbifold
$\mathcal{O}/\Stab(\mathcal{O},M)$. We can view $\mu$ either as an
{\it invariant} function on the orbit $\mathcal{O}$ or as a
function on the orbifold $\mathcal{O}/\Stab(\mathcal{O},M)$.

We now define the \emph{canonical measure} $\d\omega$.  The map
$\Im \, X(M) \to \mathcal{O}$ sends $\chi \mapsto \omega \otimes
\chi$; the map $\Im \, X(M) \to \Im \, X(A_M)$ is determined by
restriction. Let $(Y_i,\mathcal{B}_i,\mu_i)$ be finite measure
spaces with $i=1,2$ and let $f:Y_1 \to Y_2$ be a measurable map.
Then $\mu_1$ is the \emph{pull-back} of $\mu_2$ if $\mu_1(f^{-1}E)
= \mu_2(E)$ for all $E \in \mathcal{B}_2$.  This surely is the
meaning of \emph{pr\'{e}serve localement les mesures} in
\cite[p.239, 302]{Wald}.

The compact group $\Im \, X(A_M)$ is assigned the Haar measure of
total mass $1$.  Choose Haar measure on the compact orbit
$\mathcal{O}$.  Now $\Im \, X(M)$ admits two pull-back measures:
\[\Im \, X(A_M) \leftarrow \Im \, X(M) \rightarrow \mathcal{O}.\]
These must coincide: this fixes the Haar measure $\d\omega$ on
$\mathcal{O}$, see \cite[p. 239, 302]{Wald}.

Let $E$ be a Borel set in $\mathcal{O}$ which is also a
fundamental domain for the action of $\Stab(\mathcal{O},M)$ on
$\mathcal{O}$. Since $F(\omega):\,=
\mu_{G|M}(\omega)\,d(\omega)\,\theta^G_{\omega}(\lambda(g)f^{\vee})$
is $\Stab(\mathcal{O},M)$-invariant, we have
\[
|\Stab(\mathcal{O},M)|^{-1} \cdot \int_{\mathcal{O}}
F(\omega)\d\omega = \int_E F(\omega)\d\omega.
\]
The \emph{Plancherel density}, with respect to the canonical
measure $\d\omega$, is therefore \[c(G|M)^2 \cdot \gamma(G|M)^{-1}
\cdot \mu_{G|M}(\omega)\,d(\omega)\]
where $d(\omega)$ is the formal degree of $\omega$.   It is
precisely this expression which we will compute explicitly for
$\GL(n)$. To this end, we will use the following result.

\begin{thm} \label{explicitproduct} {\rm The Product Formula \cite[V.2.1]{Wald}.}
With $M = \GL(n_1) \times \cdots
\times \GL(n_k) \subset \GL(n)$ and $\omega = \omega_1 \otimes
\cdots \otimes \omega_k$ we have
\[
\mu_{G|M}(\omega) = \prod_{1 \leq j < i \leq k} \mu_{\GL(n_i +
n_j)|\GL(n_i) \times \GL(n_j)}(\omega_i \otimes \omega_j).
\]
\end{thm}

The Plancherel measure $\nu$ is determined by the equation
\[
f(g) = \int \trace\,\pi(\lambda(g)f^{\vee}) d\nu(\pi)
\]
for all $f \in \mathcal{C}(G)$.

\begin{thm} {\rm The Bernstein Decomposition \cite{P2}.}  The
Plancherel measure $\nu$ admits a canonical Bernstein
decomposition
\[
\nu = \bigsqcup \nu_{\Omega}
\]
where $\Omega$ is a component in the Bernstein variety
$\Omega(G)$.  The domain of each $\nu_{\Omega}$ is a finite union
of orbifolds of the form $\mathcal{O}/\Stab(\mathcal{O},M)$ and is
precisely a single extended quotient.
\end{thm}

We will use Theorem 2.3 to compute the Plancherel measure of the
(extended) affine Hecke algebra $\mathcal{H}(n,q)$ (see Remark~\ref{Hecke}).

\section{Calculation of the $\gamma$ factors}

\begin{thm} \label{gammaexplicite}
We have
\[\gamma(G|M)=q^{-2\sum_{1\le i<j\le k}n_in_j}\,\frac{|\GL(n,q)|}{
|\GL(n_1,q)|\times\cdots\times|\GL(n_k,q)|}.\]
\end{thm}
\begin{proof}
By applying the formula given in \cite[p.241, l.7]{Wald} to the
group $H=\I_n+\varpi\Mat(n,\integers_F)$, we obtain
\[\gamma(G|M)=q^{-2R}\;\frac{\vol(M\cap
H)}{\vol(H)},\] with $R=\Sigma(G)^+-\Sigma(M)^+$, where
$\Sigma(G)^+$ (resp. $\Sigma(M)^+$) denotes the set of positive
roots in $G$ (resp. $M$). We have
\[
R=\sum_{1\le i<j\le k}n_in_j.\] On the other hand, since the Haar
measure on $G$ is normalized so that the volume of $\cK$ is $1$,
it follows from the exact sequence
\[1\to H\to \cK\to\GL(n,q),\]
that \[\vol(H)=|\GL(n,q)|^{-1}\;\text{ and }\;\vol(H\cap M)=
|\GL(n_1,q)|^{-1}\times\cdots\times|\GL(n_k,q)|^{-1}.\]
\end{proof}

\begin{remark}
{\rm Observe that $2\sum_{1\le i<j\le k}n_in_j$ equals the length
of the element $w=w_Mw_{\GL(n)}$, where $w_M$ (resp. $w_{\GL(n)}$)
denotes the longest element in the Weyl group of $M$ (resp.
$\GL(n)$). Let $P_{S_n}(X)$ denote the Poincar\'e polynomial of
the Coxeter group $S_n$. Then, using the fact that (see for
instance \cite[(2.6)]{M3})
\begin{equation} \label{PSnqinverse}
P_{S_n}(q^{-1})=\frac{|\GL(n,q)|}{q^{n^2-n}(q-1)^n},\end{equation}
we obtain from Theorem~\ref{gammaexplicite}
\begin{equation} \label{gammaPoincare}
\gamma(G|M)=\frac{P_{S_n}(q^{-1})}{P_{S_{n_1}}(q^{-1})\times
\cdots\times P_{S_{n_k}}(q^{-1})}.
\end{equation}

This gives the following expression for the $c$-function defined
in~\cite[I.1]{Wald}:
\begin{equation} \label{cGM}
c(G|M)=\frac{\displaystyle\prod_{1\le i<j\le
k}P_{S_{n_i+n_j}}(q^{-1})} {\displaystyle
P_{S_n}(q^{-1})\cdot\prod_{i=1}^k (P_{S_{n_i}}(q^{-1}))^{k-2}}.
\end{equation}}
\end{remark}

\section{The Langlands-Shahidi formula}

Let $\varpi$ denote a fixed uniformizer. We will choose a continuous additive
character $\Psi$ such that the conductor of $\Psi$ is
$\integers_F$. Note that Shahidi uses precisely this normalization
in \cite{Sh3}.  We shall need the $L$-factor
$L(s,\pi_1\times\pi_2)$ and the $\epsilon$-factor
$\epsilon(s,\pi_1\times\pi_2,\Psi)$ for pairs, where $s$ denotes a
complex variable (see \cite{JPSS} and \cite{Sh1}). We define the
conductor $f(\pi_1\times\pi_2)$ (see \cite{BHKcond}) and the
$\gamma$-factor $\gamma(s,\pi_1\times\pi_2,\Psi)$ (see
\cite[p.~374]{JPSS}) for pairs as
\begin{equation} \label{conducpairs}
\epsilon(0,\pi_1\times\pi_2,\Psi)=q^{f(\pi_1\times\pi_2)}\, \cdot\,
\epsilon(1,\pi_1\times\pi_2,\Psi),
\end{equation}
\begin{equation} \label{gammapairs}
\gamma(s, \pi_1 \times \pi_2,\Psi) =
\epsilon(s,\pi_1 \times \pi_2,\Psi) \cdot L(1-s, \pi_1^{\vee}
\times \pi_2^{\vee})/L(s,\pi_1 \times \pi_2).
\end{equation}

We assume in this section that $P$ is the upper block triangular
maximal parabolic subgroup of $G$ with Levi subgroup
$M=\GL(n_1)\times\GL(n_2)$. We have the Langlands-Shahidi formula
for the Harish-Chandra $\mu$-function, see \cite[\S 7]{Sh2} or
\cite[\S 6]{Sh3}:
\begin{equation} \label{LSformula}
\mu_{G|M}(\omega_1 \otimes \omega_2) = \gamma(G|M)^2 \cdot
\frac{\gamma(0, \omega_1^{\vee} \times \omega_2, \Psi)}{\gamma(1,
\omega_1^{\vee} \times \omega_2, \Psi)}.
\end{equation}

It is useful to note that
\begin{equation} \label{LS2}
\frac{\gamma(0, \omega_1^{\vee} \times \omega_2, \Psi)}{\gamma(1,
\omega_1^{\vee} \times \omega_2, \Psi)}=  q^{f(\omega_1^{\vee}
\times \omega_2)} \cdot L''
\end{equation}
where
\begin{equation} \label{LS3}
L'' = \frac{L(1,\omega_1 \times \omega_2^{\vee})L(1,
\omega_1^{\vee} \times \omega_2)}{L(0,\omega_1 \times
\omega_2^{\vee})L(0, \omega_1^{\vee} \times \omega_2)}.
\end{equation}

For any smooth representation $\pi$ of $G$ and any quasicharacter
$\chi$, we denote by $\chi\pi$ the twist of $\pi$ by $\chi$:
\[\chi\pi:=(\chi\circ\det)\otimes\pi.\]
If $\sigma_1$ (resp. $\sigma_2$) is an irreducible supercuspidal
representation of $\GL(m_1)$ (resp. $\GL(m_2)$), then we have
$L(s,\sigma_1 \times \sigma_2^{\vee}) = 1$ unless $\sigma_1 \cong
\chi\sigma_2$ with $\chi$ an unramified quasicharacter of
$F^{\times}$.

\smallskip

The next formula is from \cite[p.~292]{Sh2} or \cite[Prop.~8.1]{JPSS}.
\begin{lem} \label{Lcusp}  Let
$\sigma_2$ have torsion number $r$ and let $\sigma_1 \cong
\chi\sigma_2$ with $\chi$ an unramified quasicharacter such that
$\chi(\varpi) = \zeta$. Then we have
\[
L(s, \sigma_1 \times \sigma_2^{\vee}) =(1 - \zeta^{-r} q^{-rs})^{-1}.
\]
\end{lem}

\smallskip

Let $\chi_1, \chi_2$ be unramified (unitary) characters of $F^{\times}$.
The group of unramified (unitary) characters $\Im \, X(M)$ of $M$ has, via the
map
\[
(\chi_1 \circ \det) \otimes (\chi_2 \circ \det) \mapsto
(\chi_1(\varpi), \chi_2(\varpi))\] the structure of the compact
torus $\mathbb{T}^2$.

Let $\pi_i$ be in the discrete series of $\GL(n_i)$ with $i =
1,2$, and let $\pi_i$ have torsion number $r$. Consider now the
{\it orbit} $\Im X(M) \cdot (\pi_1 \otimes \pi_2)$ in the
Harish-Chandra parameter space $\Omega^\temp(G)$. The action of
$\Im X(M)$ creates a short exact sequence
\[
1 \to \mathcal{G} \to \mathbb{T}^2 \to \mathbb{T}^2 \to 1
\]
with
\[
\mathbb{T}^2 \to \mathbb{T}^2,\; (\zeta_1, \zeta_2) \mapsto
(\zeta_1^r, \zeta_2^r).
\]
The finite group $\mathcal{G}$ is precisely the finite group in
\cite[Lemma 25]{Be} and is the product of cyclic groups:
\[
\mathcal{G} = \mathbb{Z}/r\mathbb{Z} \times
\mathbb{Z}/r\mathbb{Z}.
\]
We will write $z_1 = \zeta_1^r, z_2 = \zeta_2^r$ so that $z_1,
z_2$ are precisely the coordinates of a point in the orbit.

\begin{remark} \label{Langlandsquotients}
{\rm We recall the following facts about the discrete series of
$\GL(n)$. Let $\pi_1$ and $\pi_2$ be two discrete series
representations of $\GL(n_1)$ and $\GL(n_2)$, respectively. By
\cite{Ze}, there exist two pairs of integers $(m_1,l_1)$ and
$(m_2,l_2)$ and two irreducible unitary supercuspidal
representations $\sigma_1$ and $\sigma_2$ of $\GL(m_1)$ and
$\GL(m_2)$ respectively such that, for $i=1$, $2$, we have
$l_im_i=n_i$ and the representation $\pi_i$ is the unique
irreducible quotient associated to the Zelevinsky segment
\[
\{|\det|^{-g_i}\sigma_i,|\det|^{-g_i+1}\sigma_i,\ldots,
|\det|^{g_i-1}\sigma_i, |\det|^{g_i}\sigma_i\},\]
where $2g_i+1=l_i$.
We will follow the notation in Arthur-Clozel \cite[p.~61]{AC} and
write
\[
\pi_i = \St(\sigma_i,l_i).\] So $\pi_i$ is a \emph{generalized
Steinberg representation}.  We observe that
\[\chi\pi_i=\St(\chi\sigma_i,l_i).\]
It follows that the torsion numbers of $\pi_i$ and $\sigma_i$ are
equal.}
\end{remark}

\begin{thm} \label{mumaxs}
Let $\sigma_1$, $\sigma_2$ be irreducible unitary supercuspidal
representations of $\GL(m_1), \GL(m_2)$. Let $\pi_1$, $\pi_2$ be
discrete series representations of $\GL(n_1)$, $\GL(n_2)$ such
that $\pi_i = \St(\sigma_i,l_i)$. Let $\chi_1$,
$\chi_2$ be unramified characters. If $\sigma_1 \neq
\chi\sigma_2$ for any unramified quasicharacter $\chi$ of $F^\times$ then,
as a function on the compact torus
$\mathbb{T}^2$, $\mu_{G|M}(\chi_1 \pi_1 \otimes \chi_2 \pi_2)$ is
constant: we have
\[
\mu_{G|M}(\chi_1 \pi_1 \otimes \chi_2 \pi_2) =  \gamma(G|M)^2
\cdot q^{l_1l_2f(\sigma_1^{\vee} \times \sigma_2)}
\]
We also have
\[
f(\pi_1^{\vee} \times \pi_2) = l_1l_2f(\sigma_1^{\vee} \times
\sigma_2).
\]
\end{thm}
\begin{proof}
Let $\omega_i=\chi_i\pi_i$ and $\tau_i=\chi_i\sigma_i$ for
$i=1,2$. We will use the multiplicative property of the
$\gamma$-factors.  From \cite[Th.~3.1]{JPSS} or \cite[p.~254]{HT},
we have, with $b= g_1+g_2$,
\[\gamma(s,\omega_1^{\vee}\times \omega_2, \Psi)=
\prod_{i=0}^{l_1-1}\prod_{j=0}^{l_2-1}\gamma(s,|\;|^{i+j-b}
\tau_1^{\vee} \times \tau_2,\Psi).\]
On the other hand $\gamma(s,|\;|^{i+j-b} \tau_1^{\vee} \times
\tau_2,\Psi)$ equals
\[\epsilon(s,|\;|^{i+j-b}
\tau_1^{\vee} \times
\tau_2,\Psi)\cdot\frac{L(1-s,|\;|^{-i-j+b} \tau_1 \times
\tau_2^{\vee})}{L(s,|\;|^{i+j-b} \tau_1^{\vee} \times
\tau_2)}.\]
Since
\[\frac{\epsilon(0,|\;|^{i+j-b}
\tau_1^{\vee} \times \tau_2,\Psi)} {\epsilon(1,|\;|^{i+j-b}
\tau_1^{\vee} \times \tau_2,\Psi)}=
q^{f(|\;|^{i+j-b}\tau_1^{\vee} \times \tau_2)} =
q^{f(\tau_1^{\vee} \times \tau_2)}
=q^{f(\sigma_1^{\vee} \times \sigma_2)},
\] it follows that
\begin{equation} \label{gamma}
\frac{\gamma(0,\omega_1^{\vee}\times \omega_2, \Psi)}
{\gamma(1,\omega_1^{\vee}\times \omega_2, \Psi)}= q^{l_1l_2 \cdot
f(\sigma_1^{\vee} \times \sigma_2)}\cdot L', \end{equation} with
\begin{equation} \label{Lprime}
L'=\prod_{i=0}^{l_1-1}\prod_{j=0}^{l_2-1}
\frac{L(1,|\;|^{-i-j+b} \tau_1 \times
\tau_2^{\vee})}{L(0,|\;|^{i+j-b} \tau_1^{\vee} \times
\tau_2)}\cdot \frac{L(1,|\;|^{i+j-b} \tau_1^{\vee} \times
\tau_2)}{L(0,|\;|^{-i-j+b} \tau_1 \times \tau_2^{\vee})}.
\end{equation}
Since $\sigma_1 \neq\chi \sigma_2$, then $\tau_1 \neq\chi\tau_2$
for any unramified quasicharacter $\chi$, and $L' = 1$.

The multiplicative property of the $L$-factors \cite[Theorem
8.2]{JPSS} implies that $L'' = 1$.   Therefore, by (\ref{LS2}) we
have
\begin{equation} \label{gamma2}
\frac{\gamma(0,\omega_1^{\vee}\times \omega_2, \Psi)}
{\gamma(1,\omega_1^{\vee}\times \omega_2, \Psi)}=
q^{f(\omega_1^{\vee} \times \omega_2)}
\end{equation}

Then the results follow from the Langlands-Shahidi
formula~(\ref{LSformula}), and from ~(\ref{gamma}) and
~(\ref{gamma2}).
\end{proof}

\begin{thm} \label{mumaximalparabolic}
Let $\sigma$ be an irreducible unitary supercuspidal
representations of $\GL(m)$ with torsion number $r$.
Let $\pi_1$, $\pi_2$ be discrete series representations of $\GL(n_1)$,
$\GL(n_2)$, with $n_i=l_im$, such that $\pi_i=\St(\sigma,l_i)$.
Let $\chi_1$, $\chi_2$ be unramified characters. Let $\chi_i(\varpi)=\zeta_i$,
$z_i = \zeta^r_i$, $i=1,2$.
Then, as a function on the
compact torus $\mathbb{T}^2$ with co-ordinates $(z_1,z_2)$, we
have
\[
\mu_{G|M}(\chi_1 \pi_1 \otimes \chi_2 \pi_2) =  \gamma(G|M)^2
\cdot q^{l_1l_2f(\sigma^{\vee} \times \sigma)} \cdot \prod
\left|\frac{1 -
z_2z_1^{-1}q^{gr}}{1-z_2z_1^{-1}q^{-(g+1)r}}\right|^2
\]
where the product is over those $g$ for which $|g_1-g_2| \leq g
\leq g_1+g_2$. Note that $g_1-g_2$ and $g_1+g_2$ can both be half
integers.

We also have
\[
f(\pi_1^{\vee} \times \pi_2) = l_1l_2f(\sigma^{\vee} \times
\sigma) + r(l_1l_2 - \min(l_1,l_2)).
\]
\end{thm}
\begin{proof}
Let $\tau_i=\chi_i\sigma$. We have
\[L'=\prod_{i=0}^{l_1-1}\prod_{j=0}^{l_2-1}
\frac{L(1- i - j +b,\tau_1 \times \tau_2^{\vee})}{L(i + j-b
,\tau_1^{\vee} \times \tau_2)}\cdot \frac{L(i + j
+1-b,\tau_1^{\vee} \times \tau_2)}{L(-i - j +b, \tau_1
\times \tau_2^{\vee})},\]
where $L'$ is defined by~(\ref{Lprime}).

Now we delve into the combinatorics. To this end, we make a change
of variable, and a change of notation.

Let $\lambda(s) = L(s, \tau_1^{\vee} \times \tau_2)$,
$\lambda^*(s) = L(s, \tau_1 \times \tau_2^{\vee})$. Note that, for
all $s \in \mathbb{R}$, $\lambda^*(s)$ is the complex conjugate of
$\lambda(s)$. Let now $k = i+j-b$. We have
\[L'=\prod_{i=0}^{l_1-1}\prod_{j=0}^{l_2-1}
\frac{\lambda^*(1- k)}{\lambda(k)}\cdot
\frac{\lambda(1+k)}{\lambda^*(-k)}.\]
We now define the function \[a : \{-b, -b+1, \ldots, b-1,b\}
\longrightarrow \{1,2,3,\ldots,\min(l_1,l_2)\}\] as follows:
\[ a(k) = \sharp\{(i,j): k = i+j-b, 0 \leq i \leq l_1-1, 0\leq j
\leq l_2-1\}.\] Note that the function $a$ is even: $a(-k) =
a(k)$.   It first increases, then is constant with its maximum
value $\min(l_1,l_2)$, then decreases. Quite specifically, we have
\begin{itemize}
\item $a(-b) = 1$
\item $-b \leq k < -|g_1 - g_2| \Rightarrow a(k+1) - a(k) = 1$
\item $a(-|g_1-g_2|) = \min(l_1,l_2)$
\item $- |g_1 - g_2| \leq k < |g_1 - g_2| \Rightarrow a(k+1) =
a(k)$
\item $a(|g_1-g_2|) = \min(l_1,l_2)$
\item $|g_1 - g_2| \leq k < b \Rightarrow a(k+1) - a(k) =
-1$
\item $a(b) = 1$.
\end{itemize}
We have
\begin{eqnarray}
L' & =& \prod_{k = -b}^{b} \frac{\lambda^*(1-
k)^{a(k)}}{\lambda(k)^{a(k)}}\cdot
\frac{\lambda(1+k)^{a(k)}}{\lambda^*(-k)^{a(k)}}\nonumber\\& = &
\prod_{k = -b}^{b} \frac{\lambda^*(1 +
k)^{a(k)}}{\lambda(k)^{a(k)}}\cdot
\frac{\lambda(1+k)^{a(k)}}{\lambda^*(k)^{a(k)}}\nonumber\\& = &
\prod_{k = -b}^{b}\left| \frac{\lambda(1 +
k)^{a(k)}}{\lambda(k)^{a(k)}}\right|^2.
\end{eqnarray}
We also have, setting $a(1+b) = 0$,
\begin{eqnarray}
 \prod_{k = -b}^{b} \frac{\lambda(1+
k)^{a(k)}}{\lambda(k)^{a(k)}} & = & \frac{1}{\lambda(-b)} \cdot
\prod_{k = -b}^{b} \frac{\lambda(1+
k)^{a(k)}}{\lambda(1+k)^{a(1+k)}}\nonumber\\ & = &
\frac{1}{\lambda(-b)} \prod_{k = -b}^{-|g_1-g_2|-1}
\frac{1}{\lambda(k+1)}\cdot \prod_{k = |g_1 - g_2|}^{b}
\lambda(1+k)\nonumber\\& = & \frac{\lambda(1+b)}{\lambda(-b)}
\cdots \frac{\lambda(1 + |g_1 -
g_2|)}{\lambda(-|g_1-g_2|)}\nonumber\\& = & \prod_{g =
|g_1-g_2|}^{g_1+g_2}\frac{\lambda(1+g)}{\lambda(-g)}.
\end{eqnarray}

Note that $\tau_2 = \chi \tau_1$ where $\chi(\varpi) =
\zeta_2\zeta_1^{-1}$.  Therefore $\chi(\varpi)^{-r} =
z_1z_2^{-1}$. The first result now follows immediately from Lemma
4.1, since
\[
\lambda(g) = L(g,\tau_1^{\vee} \times \tau_2)= L(g, \tau_2 \times
\tau_1^{\vee}) = (1 - z_1z_2^{-1}q^{-gr})^{-1}.
\]
Note also that $|1 - z_2z_1^{-1}q^{-gr}| = |1 -
z_1z_2^{-1}q^{-gr}|$ since $z_2z_1^{-1}, z_1z_2^{-1}$ are complex
conjugates.

In addition we have
\[
|1-z_2z_1^{-1}q^{gr}|^2 = |q^{gr}-z_2z_1^{-1}|^2 =
q^{2gr}|1-z_2z_1^{-1}q^{-gr}|^2
\]
and so we have
\[\left|\frac{\lambda(g)}{\lambda(-g)}\right|^2 = q^{2gr}.\]

The multiplicative property of the $L$-factors \cite[Theorem
8.2]{JPSS} leads to the equation \[ L'' = \prod_{g =
|g_1-g_2|}^{g_1+g_2}\left|\frac{\lambda(1+g)}{\lambda(g)}\right|^2\]

Therefore we have
\begin{eqnarray}
L'/L'' & = & \prod_{g =
|g_1-g_2|}^{g_1+g_2}\left|\frac{\lambda(g)}{\lambda(-g)}\right|^2\nonumber\\
 &= &  \prod_{g =
|g_1-g_2|}^{g_1+g_2}q^{2rg}\nonumber\\
& = & q^{r(l_1l_2 - \min(l_1,l_2))}
\end{eqnarray} thanks to the identity
\[
2|g_1-g_2| + \cdots + 2(g_1+g_2) = l_1l_2-\min(l_1,l_2)
\]
which follows from the classic identity
\[
2|g_1-g_2| +1 + \cdots +2(g_1+g_2)+1 = l_1l_2.
\]
Since
\[ \frac{\gamma(1,\omega_1^{\vee} \times
\omega_2,\psi_F)}{\gamma(0,\omega_1^{\vee} \times
\omega_2,\psi_F)} = q^{f(\omega_1^{\vee} \times \omega_2)}\cdot
L'' = q^{l_1l_2f(\sigma^{\vee} \times \sigma)}\cdot L'
\]
we have
\[q^{f(\omega_1^{\vee} \times \omega_2)} = q^{l_1l_2f(\sigma^{\vee} \times \sigma)}\cdot L'/L''
= q^{l_1l_2f(\sigma^{\vee} \times \sigma)} q^{r(l_1l_2 -
\min(l_1,l_2))}
\]

and we conclude that
\[
f(\pi_1^{\vee} \times \pi_2) = l_1l_2f(\sigma^{\vee} \times
\sigma) + r(l_1l_2 - \min(l_1,l_2)).
\]

\end{proof}

The above formulas are invariant under the map $(z_1,z_2) \mapsto
(\lambda z_1, \lambda z_2)$ with $\lambda$ a complex number of
modulus $1$, and under the map $(z_1,z_2) \mapsto (z_2,z_1)$. In
section~6 of the paper we shall interpret $q^r$ as the cardinality
$q_K$ of the residue field of a canonical extension field $K/F$.

For example, let $M = \GL(1) \times \GL(2) \subset \GL(3)$,
$\omega_1 = 1$, $\omega_2 = \St(2)=\St(1,2)$.
We have $l_1 = 1$, $l_2 = 2$, $g_1 =
0$, $g_2 = 1/2$, $r = 1$. This gives the following (rational) function
on the $2$-torus:
\[\mu(\chi_1 \otimes \chi_2\, \St(2)) = \gamma(\GL(3)|M)^2 \cdot q \cdot
\left|\frac{1-z_2z_1^{-1}q^{-1/2}}{1-z_2z_1^{-1}q^{-3/2}}\right|^2
.\]

\begin{thm}  Let $G = \GL(2m), M = \GL(m) \times \GL(m)$ and let
$\sigma$ be an irreducible unitary supercuspidal representation of
$\GL(m)$ with torsion number $r$. Then we have \[ \mu_{G|M}(\chi_1
\sigma \otimes \chi_2 \sigma) = \gamma(G|M)^2 \cdot
q^{f(\sigma^{\vee} \times \sigma)} \cdot \left|\frac{1 -
z_2z_1^{-1}}{1-z_2z_1^{-1}q^{-r}}\right|^2
\]
\end{thm}
\begin{proof} This follows from Theorem 4.4 by taking $l_1 = l_2 =
1$, so that $g_1 = g_2 = g = 0$.
\end{proof}

\section{The Bernstein decomposition of Plancherel measure}

\subsection{The one exponent case}

Let $X$ be a space on which the finite group $\Gamma$ acts. The
extended quotient associated to this action is the
quotient space $\tilde{X}/\Gamma$ where
\[
\tilde{X} = \{(\gamma,x) \in \Gamma \times X : \gamma x = x \}.
\]
The group action on $\tilde{X}$ is $g.(\gamma, x) = (g\gamma
g^{-1}, gx)$.   Let $X^{\gamma} = \{ x \in X : \gamma x = x \}$
and let $Z(\gamma)$ be the $\Gamma$-centralizer of $\gamma$. Then
the extended quotient is given by:
\[
\tilde{X}/\Gamma = \bigsqcup_{\gamma}X^{\gamma}/Z(\gamma)
\]
where one $\gamma$ is chosen in each $\Gamma$-conjugacy class. If
$\gamma = 1$ then $X^{\gamma}/Z(\gamma) = X/\Gamma$ so {\it the
extended quotient always contains the ordinary quotient}:
\[
\tilde{X}/\Gamma = X/\Gamma \sqcup \ldots
\]
We shall need only the special case in which $X$ is the compact
torus $\mathbb{T}^n$ of dimension $n$ and $\Gamma$ is the
symmetric group $S_n$ acting on $\mathbb{T}^n$ by permuting
co-ordinates.

Let $\beta$ be a partition of $n$, and let $\gamma$ have cycle
type $\beta$.   Each cycle provides us with one circle, and
cycles of equal length provide us with a symmetric product of
circles.   For example, the extended quotient
$\widetilde{\mathbb{T}^5}/S_5$ is the following disjoint union of
compact orbifolds (one for each partition of $5$):
\[
\mathbb{T} \sqcup \mathbb{T}^2 \sqcup \mathbb{T}^2 \sqcup
(\mathbb{T} \times \Sym^2 \mathbb{T}) \sqcup (\mathbb{T} \times
\Sym^2 \mathbb{T}) \sqcup (\mathbb{T} \times \Sym^3 \mathbb{T})
\sqcup \Sym^5 \mathbb{T}
\]
where $\Sym^n \mathbb{T}$ is the $n$-fold symmetric product of the
circle $\mathbb{T}$.  This extended quotient is a model of the
arithmetically unramified tempered dual of $\GL(5)$.

Let $\Omega \subset \Omega (\GL(n))$ have one exponent $e$. Then
we have $e|n$ and so $em = n$.

There exists an irreducible unitary supercuspidal representation
$\sigma$ of $\GL(m)$ such that the conjugacy class of the cuspidal
pair $(\GL(m) \times \dots \times \GL(m), \sigma \otimes \dots
\otimes \sigma)$ is an element in $\Omega$. We have $\Omega \cong
\Sym^e \, \mathbb{C}^{\times}$ as complex affine algebraic
varieties. Consider now a partition $p=(l_1,\ldots,l_k)$ of $e$ into $k$ parts,
and write $2g_1 + 1 = l_1$, $\ldots$, $2g_k + 1 = l_k$. Let
\[
\pi_i = \St(\sigma,l_i)\] as in Remark 3.2. Then
$\pi_1 \in \cE_2(\GL(ml_1))$, $\ldots$, $\pi_k \in
\cE_2(\GL(ml_k))$. Note that $ml_1 + \ldots + ml_k = n$ so that
$\GL(ml_1) \times \ldots \times \GL(ml_k)$ is a standard Levi
subgroup $M$ of $\GL(n)$.  Now consider
\[
\pi = \chi_1\pi_1 \otimes \ldots \otimes \chi_k\pi_k
\]
with $\chi_1, \ldots, \chi_k$ unramified (unitary) characters.
Then $\pi \in \cE_2(M)$.  We have
\[
\omega = \I^{\GL(n)}_{MN}(\pi \otimes 1) \in \Irrt\, \GL(n)
\]
and each element $\omega \in \Irrt \, \GL(n)$ for which
$inf.ch.\omega \in \Omega$ is accounted for on this way.  As
explained in detail in \cite{P2}, we have
\begin{equation} \label{equiv}
\widetilde{X}/\Gamma \cong \Irrt \, \GL(n)_{\Omega}
\end{equation}
where $X = \mathbb{T}^e$, $\Gamma = S_e$, \ie
\[
\bigsqcup_{\gamma} X^{\gamma}/Z(\gamma) \cong \Irrt \,
\GL(n)_{\Omega}.
\]

The partition $p=(l_1,\ldots,l_k)$ of $e$ determines a permutation
$\gamma$ of the set $\{1,2, \ldots, e\}$: $\gamma$ is the product
of the cycles $(1,\ldots,l_1) \cdots (1, \ldots, l_k)$. Then the
fixed set $X^{\gamma}$ is
\[
\{(z_1, \ldots, z_1, \ldots, z_k, \ldots,z_k) \in \mathbb{T}^e :
z_1, \ldots, z_k \in \mathbb{T}\}
\]
and so $X^{\gamma} \cong \mathbb{T}^k$.

Explicitly, we have
\[
X^{\gamma} \longrightarrow \Irrt\,\GL(n)_{\Omega}
\]
\[
(z_1, \ldots,z_k) \mapsto \I^{\GL(n)}_{MN}(\chi_1\pi_1 \otimes
\cdots \otimes \chi_k\pi_k)
\]
with $\chi_1(\varpi) = \zeta_1, \ldots, \chi_k(\varpi) =
\zeta_k, z_1 = \zeta_1^r, \ldots, z_k = \zeta_k^r$ exactly as
in Theorem 3.2.  This map is constant on each $Z(\gamma)$-orbit
and descends to an {\it injective} map
\[
X^{\gamma}/Z(\gamma) \to \Irrt\,\GL(n)_{\Omega}.
\]

Taking one $\gamma$ in each $\Gamma$-conjugacy class we have the
bijective map
\[
\bigsqcup_{\gamma}X^{\gamma}/Z({\gamma})\cong \Irr\,
\GL(n)_{\Omega}.
\]

This bijection, by transport of structure, equips $\Irrt \,
\GL(n)_{\Omega}$ with the structure of disjoint union of finitely
many compact orbifolds.

We now describe the restriction $\mu_{\Omega}$ of Plancherel
density to the compact orbifold $X^{\gamma}/Z(\gamma)$.

\begin{thm} \label{calculationmufunction}
Let $\sigma$ be an irreducible unitary supercuspidal
representation of $\GL(m)$ with torsion number $r$. For $i=1$,
$\ldots$, $k$, let \[\pi_i = \St(\sigma,l_i),\] let
$\chi_i$ be an unramified character with $\chi_i(\varpi)=\zeta_i$,
and let $z_i=\zeta^r_i$.

Then, as a
function on the compact torus $\mathbb{T}^k$ with co-ordinates
$(z_1$, $\ldots$, $z_k)$ we have
\[
\mu(\chi_1\pi_1 \otimes \cdots \otimes \chi_k\pi_k) = const. \prod
\left|\frac{1 - z_jz_i^{-1}q^{gr}}{1-z_jz_i^{-1}q^{-(g+1)r}}\right|^2
\]
where the product is taken over those $i$, $j$, $g$ for which the
following inequalities hold:
$1 \leq i < j \leq k$, $|g_i - g_j| \leq g \leq g_i + g_j$, $2g_i
+ 1 = l_i$.
\end{thm}
\begin{proof} Apply Theorem~\ref{mumaximalparabolic} and the Harish-Chandra
product formula, Theorem~\ref{explicitproduct}.
Note that the function
\[
(z_1, \ldots,z_k) \mapsto const. \prod \left|\frac{1 -
z_jz_i^{-1}q^{gr}}{1-z_jz_i^{-1}q^{-(g+1)r}}\right|^2
\]
is a $Z(\gamma)$-invariant function on the $\gamma$-fixed set
$X^{\gamma} = \mathbb{T}^k$, and descends to a non-negative
function on the orbifold $X^{\gamma}/Z(\gamma)$:
\[
X^{\gamma}/Z(\gamma) \longrightarrow \mathbb{R}_{+}.
\]
\end{proof}

In the above theorem, the co-ordinates $z_1, \ldots, z_k$ should
be thought of as \emph{generalized Satake parameters}.  The
$k$-tuple $t = (z_1, \ldots, z_k)$ is a point in the standard
maximal torus $T$ of the unitary group $U(k, \mathbb{C})$.  In
that case, the roots of the unitary group are given by
\[
\alpha_{ij}(t) = z_i/z_j.
\]
The $\mu$-function may now be written in the more invariant form
\[
\mu(\chi_1\pi_1 \otimes \cdots \otimes \chi_k\pi_k) = const. \prod
(1 - \alpha(t)q^{gr})(1- \alpha(t)q^{-(g+1)r})^{-1}
\]
where the product is taken over all roots $\alpha = \alpha_{ij}$
of $U(k,\mathbb{C})$ and all $g$ for which the following
inequalities hold: $1 \leq i \leq k, 1 \leq j \leq k, i \neq j$,
$|g_i - g_j| \leq g \leq g_i + g_j$, $2g_i + 1 = l_i$.

\begin{thm}
\label{numericalconstant} We have the following numerical formula
for const.
\[
const.= q^{\ldee f(\sigma^\vee\times\sigma)} \,\cdot
\,\gamma(G|M)^2\cdot c(G|M)^2 ,\] where $\ldee=\sum_{1\le i<j\le
k} l_il_j$.
\end{thm}
\begin{proof}
The numerical constant is determined by Theorem~\ref{mumaximalparabolic} and
Theorem~\ref{explicitproduct}. Explicitly, for $i,j\in\{1,\ldots, k\}$,
setting
\[\gamma_{i,j}:=
\gamma(\GL(n_i+n_j)|\GL(n_i)\times\GL(n_j)),\]
for the $\gamma$-factor of the Levi subgroup $GL(n_i)\times\GL(n_j)$ of the
maximal standard parabolic subgroup in $\GL(n_i+n_j)$,
\begin{eqnarray}
const.=& q^{\sum_{1\le i<j\le k} l_il_jf(\sigma^\vee\times\sigma)}
\,\cdot\, \prod_{1\le i<j\le k}\gamma_{i,j}^2 \nonumber\\ =&
q^{\ldee f(\sigma^\vee\times\sigma)} \,\cdot \,\gamma(G|M)^2\cdot
c(G|M)^2 .\nonumber
\end{eqnarray}
\end{proof}

\begin{corollary} \label{ShahidiIII}
We have
\[j(\omega)= q^{\ldee f(\sigma^\vee\times\sigma)} \cdot
\prod
\left|\frac{1-z_jz_i^{-1}q^{-(g+1)r}}{1 - z_jz_i^{-1}q^{gr}}\right|^2
\]
\end{corollary}
\begin{proof}
This follows immediately from
Theorems~\ref{calculationmufunction}, \ref{numericalconstant} and
the fact that
\[c(G|M)^{-2}\,\gamma(G|M)^{-1}\,\mu_{G|M}(\omega)\,=
\gamma(G|M)\,j(\omega)^{-1}.\]
\end{proof}

\smallskip

Given $G=\GL(n)=\GL(n,F)$ choose $e|n$ and let $m = n/e$.
Let $\Omega$ be a Bernstein component in $\Omega(\GL(n))$ with one
exponent $e$.  The compact extended quotient attached to $\Omega$
has finitely many components, each component is a compact
orbifold. We now have enough results to write down explicitly the
component $\mu_{\Omega}$.
Let $l_1 + \cdots + l_k = e$ be a partition of $e$, let $\gamma = (1,
\ldots, l_1) \cdots (1, \ldots, l_k) \in S_e = \Gamma$, $g_1 = (l_1
- 1)/2$, $\ldots$, $g_k = (l_k - 1)/2$.  Then we have the fixed set
$X^{\gamma} = \mathbb{T}^k$.
Let $\sigma$ be an irreducible unitary supercuspidal re\-pre\-sen\-ta\-tion of
the group $\GL(m)$ and let the con\-ju\-ga\-cy class of the
cus\-pi\-dal pair $(\GL(m)^e, \sigma^{\otimes e})$ be a point in the
Bernstein component $\Omega$.
Let $r$ be the torsion number of $\sigma$ and choose a field
$K$ such that $q_K=q_F^r$.

We have~(\ref{equiv}):
\[
\Irrt\GL(n,F)_{\Omega} \cong \widetilde{X}/\Gamma .
\]
This compact Hausdorff space admits the
Harish-Chandra {\it canonical measure} $d\omega$: on each
connected component in the extended quotient
$\widetilde{X}/\Gamma$, $d\omega$ restricts to the quotient by
the centralizer $Z(\gamma)$ of the normalized Haar measure on the compact torus
$X^{\gamma}$.

Let $d\nu$ denote Plancherel measure on the tempered dual of $\GL(n,F)$.

\begin{thm} \label{PlancherelHCmeasure}
On the component $X^{\gamma}/Z(\gamma)$ of the extended quotient
$\widetilde{X}/\Gamma$ we have:
\[
d\nu(\omega) = q^{\ldee f(\sigma^\vee\times\sigma)}
\,\cdot\gamma(G|M)\,\cdot\, d(\omega)
\,\cdot\,\prod
\left|\frac{1 - z_jz_i^{-1}q^{gr}}{1-z_jz_i^{-1}q^{-(g+1)r}}\right|^2 \,\cdot\,d\omega.
\]
\end{thm}
\begin{proof}
By~(\ref{HCformula}), the Plancherel measure on
$\Irrt\GL(n,F)_{\Omega}$ is given by
\[
d\nu(\omega) = c(G|M)^{-2} \,\gamma(G|M)^{-1}\, \mu(\omega)\, d(\omega)\, d\omega
\]
Then, the result follows from Theorem~\ref{calculationmufunction} and
Theorem~\ref{numericalconstant}.
\end{proof}

\smallskip

Let $T$ be the diagonal subgroup of $G$ and take for $\Omega$ the Bernstein
component in $\Omega(G)$ which contains the cuspidal pair $(T,1)$.
Then $\Omega$ has the single exponent $n$ and
parametrizes those irreducible smooth representations of $\GL(n,F)$ which admit
nonzero Iwahori fixed vectors.

\smallskip
Now let $l_1+\cdots+l_k$ be a partition of $n$, and let
\[M=\GL(l_1,F)\times\cdots\times\GL(l_k,F)\,\subset\,\GL(n,F).\]

The formal degree of the Steinberg representation $\St(l_i)$ is
given by
\begin{equation} \label{Steinberg}
d(\St(l_i))=\frac{q^{(l_i-l_i^2)/2}}{l_i}\cdot\frac{|\GL(l_i,q)|}{q^{l_i}-1}
=\frac{1}{l_i}\cdot \prod_{j=1}^{l_i-1}(q^j-1)\end{equation}

We also have the inner product identity in pre-Hilbert space:
\[
\langle(\sigma_1 \otimes \cdots \otimes  \sigma_k)(g)\xi_1  \otimes
\cdots \xi_k, \xi_1 \otimes \cdots \otimes \xi_k\rangle \; = \;\prod
\langle\sigma_j(g)\xi_j, \xi_j\rangle.
\]
Let each $\xi_j \in V_j$ be a unit vector. With respect to the
standard normalization of all Haar measures we then have (cf.
\cite[(7.7.9)]{BK})
\[1/d_{\sigma_1 \otimes \cdots \otimes \sigma_k} = \prod \int |<\sigma_j(g)\xi_j, \xi_j>|^2
d \dot{\vol}_j= \prod 1/d_{\sigma_j}\] and so
\begin{equation} \label{productfd}
d_{\sigma_1 \otimes \cdots \otimes \sigma_k} = \prod
d_{\sigma_j}.
\end{equation}

Using (\ref{productfd}) and Theorem~\ref{gammaPoincare}, we obtain the
following result.

\begin{corollary} \label{Iwahori}
On the orbifold $X^\gamma/Z(\gamma)$ we have
\[
d\nu(\omega) = \gamma(G|M)\,\cdot\, d(\omega)\,\cdot\,
\prod
\left|\frac{1 - z_jz_i^{-1}q^{g}}{1-z_jz_i^{-1}q^{-(g+1)}}\right|^2 \,\cdot\,d\omega
\]
where \[d(\omega) = \prod d(\St(l_i)).\]

So we have
\begin{eqnarray}
d\nu(\omega) &= \displaystyle\gamma(G|M)\,\cdot\,
\prod_{i=1}^k\frac{1}{l_i}\prod_{j=1}^{l_i-1}( q^j-1)
\,\cdot\,\prod
\left|\frac{1 - z_jz_i^{-1}q^{g}}{1-z_jz_i^{-1}q^{-(g+1)}}\right|^2
\,\cdot\,d\omega\cr
&=\displaystyle\prod_{i=1}^k\frac{q^{\frac{l_i^2-l_i}{2}}(q-1)^{l_i}}{l_i(q^{l_i}-1)}\,
\cdot\,P_{S_n}(q^{-1})\,\cdot\,\prod
\left|\frac{1 - z_jz_i^{-1}q^{g}}{1-z_jz_i^{-1}q^{-(g+1)}}\right|^2
\,\cdot\,d\omega.
\end{eqnarray}
\end{corollary}

\begin{remark} \label{Hecke}
{\rm Using \cite[Theorem~3.3]{BHK}, we obtain that the Plancherel
measure of the (extended) affine Hecke algebra ${\mathcal H}(n,q)$
is given on $X^\gamma/Z(\gamma)$ by
\[\vol(I)\,\cdot\,
\gamma(G|M)\,\cdot\,d(\omega)\,\cdot\,
\prod
\left|\frac{1 - z_jz_i^{-1}q^{g}}{1-z_jz_i^{-1}q^{-(g+1)}}\right|^2 \,\cdot\,d\omega.\]

Concerning the volume $\vol(I)$: by \cite[5.4.3]{BK} we have
\[
\vol(\GL(n,\integers_F)) = \sum_{w \in W_0} \vol(IwI) = \sum_{w \in W_0}
\vol(I)\cdot q^{\ell(w)} = P_{S_n}(q)\cdot \vol(I).\]
The explicit formula is then (using (\ref{PSnqinverse})):
\[
d\nu_{{\mathcal H}(n,q)}(\omega)
=\displaystyle\prod_{i=1}^k
\frac{q^{\frac{l_i^2-l_i}{2}}(q-1)^{l_i}}{l_i(q^{l_i}-1)}\,
\cdot\,q^{\frac{n-n^2}{2}}\,\cdot\,\prod
\left|\frac{1 - z_jz_i^{-1}q^{g}}{1-z_jz_i^{-1}q^{-(g+1)}}\right|^2
\,\cdot\,d\omega,
\]
where the second product is taken over those $i$, $j$, $g$ for which the
following inequalities hold:
$1 \leq i < j \leq k$, $|g_i - g_j| \leq g \leq g_i + g_j$, $2g_i
+ 1 = l_i$.
Note that Plancherel measure for Iwahori Hecke algebras has been
already calculated by Opdam (see \cite[2.8.3]{O}).}
\end{remark}

\smallskip
We will now consider a special case.
The $p$-adic gamma function attached to the local field $K$ (see
\cite[p. 51]{Tai}) is the following
meromorphic function of a single complex variable:
\[
\Gamma_1(\zeta) = \frac{1 - q_K^{\zeta}/q_K}{1 - q_K^{-\zeta}}.\]
We will change the variable via $s = q_K^{\zeta}$ and write
\[
\Gamma_{K}(s) = \frac{1-s/q_K}{1-s^{-1}},
\]
a rational function of $s$.  Let $s \in i\mathbb{R}$ so that $s$
has modulus $1$.  Then we have
\[
1/|\Gamma_{K}(s)|^2 = \left|\frac{1 -  s}{1 -
q_K^{-1}s}\right|^2.\]

Let $T$ be the standard maximal torus in $\GL(n)$ and let
$\widehat{T}$ denote the unitary dual of $T$.   Then $\widehat{T}$
has the structure of a compact torus $\mathbb{T}^n$ (the space of
Satake parameters) and the unramified unitary principal series of
$\GL(n)$ is parametrized by the quotient $\mathbb{T}^n/S_n$.   Let
now $t = (z_1, \dots, z_n) \in \mathbb{T}^n$.  Applying the above
formulas the Plancherel density $\mu_{G|T}$ is given by
\begin{eqnarray}
\mu_{G|T} & = & const \cdot \prod_{i<j} \left | \frac{1 -
z_jz_i^{-1}}{1 - z_jz_i^{-1}/q} \right|^2 \\ & = & const \cdot
\prod_{0<\alpha} \left | \frac{1 - \alpha(t)}{1 - \alpha(t)/q}
\right|^2 \\ & = & const \cdot \prod_{\alpha}1/\Gamma(\alpha(t))
\end{eqnarray}
where $\alpha$ is a root of the Langlands dual group $\GL(n, \mathbb{C})$
so that $\alpha_{ij}(t) =
z_i/z_j$.

For $\GL(n)$, one connected component in the tempered dual is the
compact orbifold $\mathbb{T}^n/S_n$, the symmetric product of $n$ circles.
On this component we have
the Macdonald formula \cite{M1}:
\[
d\mu(\omega_{\lambda}) = const. \cdot d\lambda /
\prod_{\alpha}\Gamma(i\lambda(\alpha^\vee))
\]
the product over all roots $\alpha$ where $\alpha^\vee$ is the
coroot.   This formula is a very special case of our formula for
$\GL(n)$.
\subsection{General case}
We now pass to the general case of a component $\Omega \subset
\Omega (\GL(n))$ with exponents $e_1$, $\ldots$, $e_t$. We first
note that each component $\Omega\subset \Omega (\GL(n))$ yields up
its fundamental invariants:
\begin{itemize}
\item the cardinality $q$ of the residue field of $F$
\item
the sizes $m_i$ of the small general linear groups
\item
the exponents $e_i$
\item
the torsion numbers $r_i$
\item
the formal degrees $d_i$
\item the conductors $f_{ij} = f(\sigma_i^{\vee} \times \sigma_j)$
\end{itemize}
with $1 \leq i \leq t$.




We now construct the disjoint union
\[
E = \Omega(\GL(\infty)) = \{\bigsqcup \Omega (\GL(n)): n =
0,1,2,3, \ldots \}
\]
with the convention that $\Omega(\GL(0)) = \mathbb{C}$.

We will say that two components $\Omega_1, \Omega_2 \in E$ are
{\it disjoint} if none of the irreducible supercuspidals which
occur in $\Omega_1$ is equivalent (after unramified twist) to any
of the supercuspidals which occur in $\Omega_2$. We now define a
law of composition on {\it disjoint components} in $E$.  With the
cuspidal pair $(M_1, \sigma_1) \in \Omega_1$ and the cuspidal pair
$(M_2, \sigma_2) \in \Omega_2$ we define $\Omega_1 \times
\Omega_2$ as the unique component determined by
\[ (M_1 \times M_2, \sigma_1 \otimes \sigma_2). \]

The set $E$ admits a law of composition {\it not everywhere
defined} such that $E$ is unital, commutative and associative.
Rather surprisingly, $E$ admits prime elements: the prime elements
are precisely the components with a single exponent. Each element
in $E$ admits a unique factorization into prime elements:
\[
\Omega = \Omega_1 \times \cdots \times \Omega_t.
\]

Plancherel measure respects the unique factorization into prime
elements, modulo constants.   Quite specifically, we have

\begin{thm} Let $\Omega$ have the unique factorization
\[\Omega = \Omega_1 \times \cdots \times \Omega_t
\]
so that $\Omega$ has exponents $e_1$, $\ldots$, $e_t$ and $\Omega_1$,
$\ldots$, $\Omega_t$ are pairwise disjoint prime elements with the
individual exponents $e_1$, $\ldots$, $e_t$. Let
\[
\nu = \bigsqcup \nu_{\Omega} \] denote the Bernstein decomposition
of Plancherel measure.  Then we have
\[
\nu_{\Omega} = const.\, \nu_{\Omega_1} \cdots \nu_{\Omega_t}
\]
where $\nu_{\Omega_1}$, $\ldots$, $\nu_{\Omega_t}$ are given by
Theorem~\ref{calculationmufunction} and the constant is given, in
terms of the fundamental invariants, by
Theorem~\ref{numericalconstant}.
\end{thm}

\begin{proof} In the Harish-Chandra product formula, all the
cross-terms are constant, by Theorem 4.3.
\end{proof}


\section{Transfer-of-measure, conductor, and the formal degree formulas}

\subsection{Torsion number} \label{Tn}

The theory of types of \cite{BK} produces a canonical extension
$K$ of $F$ such that $q_K=q^r$. Indeed, let $\sigma$ be an
irreducible supercuspidal representation of $\GL(m)$, and let
$(J,\lambda)$ be a maximal simple type occurring in it. Let $\Afr$
be the hereditary $\integers_F$-order in $A=M(m,F)$
and let $E=F[\beta]$ be the field extension of $F$ attached to the stratum
(see \cite[Definition 5.5.10~(iii)]{BK}). It is proved in
\cite[Lemma 6.2.5]{BK} that
\begin{equation} \label{torsionnumber}
r=\frac{m}{e(E|F)},
\end{equation} where $e(E|F)$ denotes the ramification
index of $E$ with respect to $F$. Let $B$ denote the centraliser
of $E$ in $A$. We set $\Bfr:=\Afr\cap B$. Then $\Bfr$ is a maximal
hereditary order in $B$, see \cite[Theorem 6.2.1]{BK}. Let $K$ be
an unramified extension of $E$ which normalises $\Bfr$ and is
maximal with respect to that property, as in \cite[Proposition
5.5.14]{BK}. Then $[K:F]=m$, and (\ref{torsionnumber}) gives that
$r$ is equal to the residue index $f(K|F)$ of $K$ with respect to
$F$. Thus $Q=q^r$ is equal to the order $q_K$ of the residue field
of $K$.

Also the number $Q$ is the one which occurs for the Hecke algebra
$\mathcal{H}(\GL(m),\lambda)$ associated to $(J,\lambda)$, see
\cite[Theorem 5.6.6]{BK}. Indeed,  since the order of the residue
field of $E$ is equal to $q^{f(E|F)}$, that number is
$(q^{f(E|F)})^f$, with $$f=\frac{m}{[E:F]\,e(\Bfr)},$$ where
$e(\Bfr)$ denotes the period of a lattice chain attached to $\Bfr$
as in \cite[(1.1)]{BK}. Since $\sigma$ is supercuspidal,
$e(\Bfr)=1$ (see \cite[Corollary 6.2.3]{BK}). It follows that
\begin{equation} \label{tors}
f\,\cdot\,f(E|F)=\frac{m\,\cdot
\,f(E|F)}{[E:F]}=\frac{m}{e(E|F)}=r.
\end{equation}

\subsection{Normalization of measures}
\label{normeasures}

We will relate our normalization of measures to the measures used in
\cite[(7.7)]{BK}.
Bushnell and Kutzko work with a quotient measure $\dot{\mu}$, the quotient of
$\vol_G$ by $\vol_Z$.

Let $Z$ denote the centre of $\GL(n)$.  The second isomorphism
theorem in group theory gives:
\[
JZ/Z \cong J/J \cap Z.\]

We have \[ J\cap Z= \integers_F^{\times}.\]
One way to see this would be: $J$ contains $\Afr^\times\cap B$, where
$B$ is the centralizer in $M(n,F)$ of the extension $E$. Now certainly $Z$ is
contained in $B$. On the other hand, $\Afr$ is an $\integers_F$-order so $\Afr$
certainly contains $\integers_F$. Thanks to Shaun Stevens for this remark.

Then we have
\[
JZ/Z \cong J/\integers_F^{\times}.\]

Now $J$ is a principal $\integers_F^{\times}$-bundle over
$J/\integers_F^{\times}$. Each fibre over the base $J/
\integers_F^{\times}$ has volume $1$. The quotient measure of the
base space is then given by
\begin{equation} \label{volJ}
\dot{\mu}(JZ/Z) = \vol(J).
\end{equation}

Similar normalizations are done with $G_0 = \GL(e,K)$. We also
need the corresponding quotient measure $\ddot{\mu}$ (see
\cite[(7.7.8)]{BK}). We have
\[
\ddot{\mu}(IK^{\times}/K^{\times}) = \vol_{G_0}(I).\]

Let $M= \prod \GL(n_j)$.  We have $Z_M = \prod Z_j$, $\cK = \prod \cK_j$,
with $Z_j=Z_{\GL(n_j,F)}$ and $\cK_j=\GL(n_j,\integers_F)$.
With respect to the standard normalization of all Haar measures,
we have $\vol_M = \prod \vol_j$ (where $\vol_j$ denotes
$\vol_{\GL(n_i,F)}$) and $\vol_{Z_M} = \prod \vol_{Z_j}$. This then
guarantees that
\begin{equation} \label{volM}
\dot{\vol}_M = \prod \dot{\vol}_j.
\end{equation}

\subsection{Conductor formulas (the supercuspidal case)}
We will first recall results from \cite{BHKcond} in a suitable way
for our purpose.

Let $(J^\simple,\lambda^\simple)$ be a simple type in $\GL(2m)$
with associated maximal simple type $(J,\lambda)$ (in the
terminology of \cite[(7.2.18)~(iii)]{BK}). When $(J,\lambda)$ is
of positive level, we set $J_P=(J^\simple\cap
P)H^1(\beta,\Afr)\subset J^\simple$ (in notation
\cite[(3.1.4)]{BK}), where $P$ denotes the upper-triangular
parabolic subgroup of $\GL(2m)$ with Levi component
$M=\GL(m)\times\GL(m)$, and unipotent radical denoted by $N$.
Following \cite[Theorem 7.2.17]{BK}, we define $\lambda_P$ as the
natural representation of $J_P$ on the space of $(J\cap N)$-fixed
vectors in $\lambda^\simple$. The representation $\lambda_P$ is
irreducible and
$\lambda_P\simeq\cInd_{J_P}^{J^\simple}(\lambda^\simple)$.

The pair $(J\times J,\lambda\otimes\lambda)$ is a type in $M$
which occurs in $\sigma\otimes\sigma$, and, as shown in
\cite[prop.~1.4]{BKss}, $(J_P,\lambda_P)$ is a $\GL(2m)$-cover of
$(J\times J,\lambda\otimes\lambda)$.

\begin{thm} \label{condformulae} {\rm Conductor formulas, \cite{BHKcond}.}
Let $G_0=\GL(2,K)$, let $N_0$ denote the unipotent radical of the
standard Borel subgroup of $G_0$ and let $I$ denote the standard
Iwahori subgroup of $G_0$. We will denote by $\vol_0$ the Haar measure on
$G_0$ normalized as in subsection~\ref{normeasures}.

Let $(J^{\GL(2m)},\lambda^G)$ be any $\GL(2m)$-cover of $(J\times
J,\lambda\otimes \lambda)$.

Then
\[\frac{\vol(J^G\cap N)\cdot \vol(J^G\cap\overline N)}
{\vol_0(I\cap N_0)\cdot \vol_0(I\cap\overline N_0)}=
q^{-f(\sigma^\vee\times\sigma)}
=
\frac{j(\sigma\otimes\sigma)}{j_0(1)} ,\] where $j$, $j_0$ denote
the $j$-functions for the group $G$, $G_0$ respectively.
\end{thm}
\begin{proof}
The first equality is \cite[Theorem in \S 5.4]{BHKcond}, using the
fact that $\vol_0(I\cap N_0)\cdot \vol_0(I\cap\overline N_0)=q_K^{-1}$.
The second equality is \cite[Theorem in \S 5.4]{BHKcond} (note
that in {\it loc. cit.} the normalisations haven been taken so
that $\vol(J^G\cap N)\cdot \vol(J^G\cap\overline N)=\vol_0(I\cap N_0)\cdot
\vol_0(I\cap\overline N_0)$). It also follows directly from our
Corollary~\ref{ShahidiIII}.
\end{proof}

We will now extend the above Theorem to the case of
$M=\GL(m)^{\times e}$, with $e$ arbitrary.

\begin{corollary} \label{extensiontoeblocks}
Let $M=\GL(m)^{\times e}$ with $n=em$, et $G_0=\GL(e,K)$, let
$N_0$ denote the unipotent radical of the standard Borel subgroup
of $G_0$ and let $I$ denote the standard Iwahori subgroup of
$G_0$.

Let $(J^G,\lambda^G)$ be a cover in $G=\GL(n)$ of $(J^{\times
e},\lambda^{\otimes e})$ (the existence of which is guaranteed by
\cite{BKss}).

Then
\[\frac{\vol(J^G\cap N)\cdot \vol(J^G\cap\overline N)}
{\vol_0(I\cap N_0)\cdot \vol_0(I\cap\overline N_0)}=
q^{-\frac{e(e-1)}{2}f(\sigma^\vee\times\sigma)}
=
\frac{j(\sigma^{\otimes e})}{j_0(1)} .\]
\end{corollary}
\begin{proof}
Let $M'$ be a Levi subgroup of a parabolic subgroup in $G$ such
that $P$ is a maximal parabolic subgroup of $M'$. Then,
$M'/M\simeq \GL(2m)/\GL(m)\times\GL(m)$ and
\[\vol(J^G\cap M'\cap
N)=\vol(J^{\GL(2m)}\cap\GL(2m)\cap N).\] It follows from
\cite[Proposition 8.5~(ii)]{BKtyp} that $(J^G\cap
M',\lambda^G|J^G\cap M')$ is an $M'$-cover of $(J^{\times
e},\lambda^{\otimes e})$.

Because of the unipotency of $N$, we have
\begin{equation} \label{vJN}
\vol(J^G\cap N)=(\vol(J^{\GL(2m)}\cap\GL(2m)\cap N))^{\frac{e(e-1)}{2}},
\end{equation}
and similar equalities for the three others terms. Since
$\GL(2m)\cap N$ is the unipotent radical of the parabolic subgroup
of $\GL(2m)$ with Levi $\GL(m)\times\GL(m)$, the first equality in
the Corollary follows from Theorem~\ref{condformulae}.

The second equality follows from our Corollary~\ref{ShahidiIII}.
It is also a direct consequence of Theorem~\ref{condformulae},
using the product formula for $j$ and for $j_0$ from
\cite[IV.3.~(5)]{Wald}.
\end{proof}

\subsection{Formal degree formulas}

Using Corollary~\ref{extensiontoeblocks},
we will deduce from \cite[(7.7.11)]{BK} a formula relating the formal
degree of any discrete series of $\GL(n)$ and the formal degree of a
supercuspidal representation in its inertial support.

Given $G=\GL(n)=\GL(n,F)$ choose $e|n$ and let $m=n/e$.
Let $\sigma$ be an irreducible unitary supercuspidal
representation of $\GL(m)$
and let $(J,\lambda)$ be a maximal simple type occurring in it.
Let $g=(e-1)/2$. We consider the standard Levi subgroup
$M=\GL(m)^{\times e}$ of $\GL(n,F)$ and the supercuspidal
representation
\[\sigma_M=|\det(\;)|^{-g}\sigma\otimes\cdots\otimes|\det(\;)|^{g}\sigma\]
of it.
Then $(J_M,\lambda_M)=(J^{\times e},\lambda^{\otimes e})$ is a type in $M$
occuring in $\sigma_M$.

Let $\pi = \St(\sigma,e)$
and let $(J^\simple,\lambda^\simple)$ be a simple type in $\GL(n)$ occuring in
$\pi$ (it has associated maximal simple type $(J,\lambda)$).

The following result is rather intricate, but note that only the
{\it fundamental invariants} $m, e, r, d, f(\sigma^{\vee} \times
\sigma)$ occur in it, in line with our general philosophy.

\begin{thm} \label{formaldegrees}
We have
\[\frac{d(\pi)}{d(\sigma)^e}=\frac{m^{e-1}}{r^{e-1}e}\,\cdot\,
q^{\frac{e^2-e}{2}(f(\sigma^\vee\times\sigma)+r-2m^2)}
\,\cdot\,
\frac{(q^r-1)^e}{q^{er}-1}
\,\cdot\,\frac{|\GL(em,q)|}{|\GL(m,q)|^e}.\]
\end{thm}
\begin{remark} \label{quotient}
{\rm The right-hand side in the above equality can be rewritten, by
using~(\ref{Steinberg}), as
\[
r^{1-e}\,\cdot\,\frac{(q^{em}-1)(q^r-1)^e}{(q^m-1)^e(q^{er}-1)}\,\cdot\,
q^{\frac{e^2-e}{2}(f(\sigma^\vee\times\sigma)+r-m^2)}
\,\cdot\,
\frac{\deg(\St(em))}{(\deg(\St(m)))^e}.\]
}
\end{remark}
\begin{proof}
Let $T$ denote the diagonal torus in $\GL(e,K)$ and let $I$ denote
the Iwahori subgroup of $G_0=\GL(e,K)$ attached to the Bernstein
component in $\Omega(\GL(e,K))$ which contains the cuspidal pair
$(T,1)$. Note that $I\cap T=\GL(1,\integers_K)^{\times e}$. From
\cite[(7.7.11)]{BK}, applied to the representations $\pi$ and
$\sigma$, we have
\begin{equation} \label{fdomega}
d(\pi)=\frac{\vol_0(I)}{\vol(J^\simple)}\,\cdot\,
\frac{\dim(\lambda^\simple)}{e(E|F)}\,\cdot\,d(\pi)_0,
\end{equation}
where $d(\pi)_0$ denotes the formal degree of $\pi\in\cE_2(G_0)$,
and
\begin{equation} \label{fdsigma}
d(\sigma)=
\frac{\vol(\GL(1,\integers_K))}{\vol(J)}\,\cdot\,
\frac{\dim(\lambda)}{e(E|F)}.
\end{equation}
Using (\ref{fdomega}), (\ref{fdsigma}) and (\ref{tors}), we obtain
\begin{equation} \label{ratiofd}
\frac{d(\pi)}{d(\sigma)^e}=\frac{m^{e-1}}{r^{e-1}}\,\cdot\,
\frac{\vol_0(I)}{\vol(J^\simple)}\,\cdot\,
\frac{\vol(J^{\times e})}{\vol(\GL(1,\integers_K)^{\times e})}\,\cdot\,
\frac{\dim(\lambda^\simple)}{(\dim(\lambda))^e}\,\cdot\,d(\pi)_0.
\end{equation}
We set $J_P=(J^\simple\cap P)H^1(\beta,\Afr)\subset J^\simple$, where $P$ is
the upper-triangular parabolic subgroup of $G$ with Levi component $M$, and
unipotent radical
$N$. We define $\lambda_P$ as the
natural representation of $J_P$ on the space of $(J\cap N)$-fixed vectors in
$\lambda^\simple$. The representation $\lambda_P$ is irreducible and
$\lambda_P\simeq\cInd_{J_P}^{J^\simple}(\lambda^\simple)$.
Then $(J_P,\lambda_P)$ is a $G$-cover of $(J_M,\lambda_M)$.
In the case where $(J,\lambda)$ is of zero level, we denote by
$(J^\simple,\lambda^\simple)=(J_P,\lambda_P)$ an arbitrary $G$-cover of
$(J_M,\lambda_M)$.

Since $J^\simple\cap M=J^{\times e}=J_M=J_P\cap M$,
and \[\dim(\lambda)^e=\dim(\lambda_M)=\dim(\lambda_P)=
[J^\simple:J_P]^{-1}\,\dim(\lambda^\simple),\]
(\ref{ratiofd}) gives
\[\frac{d(\pi)}{d(\sigma)^e}=\frac{m^{e-1}}{r^{e-1}}\,\cdot\,
\frac{\vol_0(I)}{\vol(J_P)}\,\cdot\,
\frac{\vol(J_M)}{\vol_0(I\cap T)}\,\cdot\,d(\pi)_0.\]
On the other hand, by applying the formula \cite[p.241, l.7]{Wald} to the group
$J$, we obtain
\begin{equation} \label{gammaWalds}
\gamma(G|M)=\frac{\vol(J_P\cap N)\cdot \vol(J_P\cap M)\cdot
\vol(J_P\cap\overline N)}{\vol(J_P)}.
\end{equation}
Similarly we have
\[\gamma(G_0|T)=\frac{\vol_0(I\cap N_0)\cdot \vol_0(I\cap T)\cdot
\vol_0(I\cap\overline N_0)}{\vol_0(I)}.\]
We then obtain
\[\frac{d(\pi)}{d(\sigma)^e}=\frac{\gamma(G|M)}{\gamma(G_0|T)}\,\cdot\,
\frac{\vol_0(I\cap N_0)\cdot \vol_0(I\cap\overline N_0)}
{\vol(J_P\cap N)\cdot \vol(J_P\cap\overline
N)}\,\cdot\,d(\pi)_0.\]
Applying Corollary~\ref{extensiontoeblocks}, we get
\[\frac{d(\pi)}{d(\sigma)^e}=\frac{m^{e-1}}{r^{e-1}}\,\cdot\,
q^{\frac{e(e-1)}{2} f(\sigma^\vee\times\sigma)}\,\cdot\,
\frac{\gamma(G|M)}{\gamma(G_0|T)}\,\cdot\,
d(\pi)_0.\]
Since Haar measure on $\GL(e,K)$ has been normalised so that the volume
of $\GL(e,\integers_K)$ is equal to one, the formal degree
of the Steinberg representation of $\GL(e,K)$ is
given as in (\ref{Steinberg}) by
\[
d(\pi)_0=\frac{q_K^{(e-e^2)/2}}{e}\cdot\frac{|\GL(e,q_K)|}{q_K^{e}-1}
.\]
On the other hand, Theorem~\ref{gammaexplicite} gives
\[\gamma(G|M)=q^{mn-n^2}\,\cdot\,\frac{|\GL(n,q)|}{|\GL(m,q)|^e}
\;\text{ and }\;\;
\gamma(G_0|T)=q^{e-e^2}\,\cdot\,\frac{|\GL(e,q_K)|}{(q_K-1)^e}.\]
The result follows.
\end{proof}

\smallskip

We will now recall the explicit formulas for $d(\pi)$ and $d(\sigma)$
from \cite{SZ}, using also \cite{Zi}. We would like to
thank Wilhelm Zink for explaining these works to us.

Let $\eta$ be the Heisenberg representation of $J^1(\beta,\Afr)$ attached
to a maximal simple type $(J(\beta,\Afr),\lambda)$ occuring in the
supercuspidal representation $\sigma$ of $\GL(m)$
(see \cite[(5.1.1), (5.5.10)]{BK}).
Let $\Pfr$ denote the Jacobson radical of $\Afr$ and let $U^i(\Afr)=
1+\Pfr^i$. Let $\pi_\beta^1$ be the compactly induced representation
$\cInd_{J^1(\beta,\Afr)}^{U^1(\Afr)}(\eta)$.
Then $\pi_\beta^1$ is irreducible, see \cite[(5.2.3)]{BK}.
More generally the restriction of $\eta$ to
$J^i(\beta,\Afr)=J^1(\beta,\Afr)\cap (1+ \Pfr^i)$
is a multiple of an irreducible representation
$\eta^i$ which induces irreducibly to a representation $\pi_\beta^i$  of
$U^i(\Afr)$ (see \cite[2.2]{Zi}).
Let $E_{-i}$ be any field such that
\[U^1(\Afr)\cdot I_{\GL(m)}(\pi_\beta^{i+1})\cdot U^1(\Afr)\, =\,
U^1(\Afr)\cdot\GL(m/[E_{-i}:F],E_{-i})\cdot U^1(\Afr),\]
where $I_{\GL(m)}(\pi_\beta^{i+1})$ denotes the intertwining of
$\pi_\beta^{i+1}$ in $\GL(m,F)$.
In particular, we have $E_0=E$.

\begin{thm} \label{explicitdegrees} {\rm Explicit formal degrees formulas,
\cite{SZ}, \cite{Zi}.}
The formal degrees of $\sigma$ and $\pi$ are respectively given by
\[
d(\sigma)\,=\,r\,\cdot\, \frac{q^m-1}{q^r-1}\cdot
q^{(r-m+\delta)/2}\cdot \deg(\St(m)),
\]
\[
d(\pi)\,=\,r\,\cdot\, \frac{q^{em}-1}{q^{er}-1}\cdot
q^{(er-em+e^2\delta)/2}\cdot \deg(\St(em)),
\]
where
\[
\delta=rm\,\cdot\,\sum_{i\ge 0} (1 - [E_{-i}:F]^{-1}).\]
\end{thm}
\begin{proof} It follows directly from \cite[Theorem
1.1]{SZ} and \cite[Corollary~6.7]{Zi}, using the fact that $r=f(K|F)$ and
$m/e(E|F)=r$.
\end{proof}

As immediate consequences, we obtain the following results.

\begin{corollary} \label{quotientdeux}
\[\frac{d(\pi)}{d(\sigma)^{e^2}}=
r^{1-e^2}\,\cdot\,\frac{(q^{em}-1)(q^r-1)^{e^2}}{(q^{er}-1)(q^m-1)^{e^2}}\,
\cdot\,q^{(e^2-e)(m-r)/2}\,\cdot\,\frac{\deg(\St(em))}{(\deg(\St(m)))^{e^2}}.\]
\end{corollary}

\begin{remark}
{\rm
We observe that the above formula extends to the general case the formula
obtained in \cite[Theorem~4.6]{CMS} in the case where $(n,p)=1$ and $F$ has
characteristic zero. The existence of such a formula was expected in
\cite[Remark~4.7]{CMS}. Our formula also extends
\cite[Theorem~VII.3.2]{Wtypes}.}
\end{remark}

\begin{corollary} \label{quotientbis}
\[\frac{d(\pi)}{d(\sigma)^e}=
r^{1-e}\,\cdot\,\frac{(q^{em}-1)(q^r-1)^e}{(q^{er}-1)(q^m-1)^e}\,
\cdot\,q^{(e^2-e)\delta/2}\,\cdot\,\frac{\deg(\St(em))}{(\deg(\St(m)))^e}.\]
\end{corollary}

The comparison of Corollary~\ref{quotientbis} with
Remark~\ref{quotient} gives the following expression for the
conductor for pairs $f(\sigma^\vee\times\sigma)$.

\begin{thm} \label{artin}
We have
\[f(\sigma^\vee\times\sigma)=\delta+m^2-r.\]
\end{thm}

\begin{remark} \label{compBHK}
{\rm In \cite[\S 6.4]{BHK} (see also \cite[6.13]{BHK}) is introduced a certain
discrimant function $C(\beta)$ and an integer $\fc(\beta)$ such that
$C(\beta)=q^{\fc(\beta)}$.
It follows from our Theorem~\ref{calculationmufunction} and
\cite[Theorem~6.5~(i)]{BHK} that
\[\fc(\beta)=\frac{[E:F]^2}{m^2}\,\cdot\,\delta.\]
}
\end{remark}

\subsection{Conductor formulas (the discrete series case)}
\label{dc}

Let $\sigma$ be an irreducible supercuspidal
representation of $\GL(m)$, and let $(J,\lambda)$ be a maximal simple type
occurring in it. Let $e|n$, and
let $l_1 + \cdots + l_k = e$ be a partition of $e$. It
determines the standard Levi subgroup
\begin{equation} \label{LeviM}
M=\GL(l_1m)\times\cdots\times\GL(l_km)\subset\GL(n,F).
\end{equation}
Let $g_1 = (l_1-1)/2$, $\ldots$, $g_k=(l_k-1)/2$, and let $\pi_1$,
$\ldots$, $\pi_k$ be discrete series representations of
$\GL(l_1m)$, $\ldots$, $\GL(l_km)$ such that
$\pi_i=\St(\sigma,l_i)$. Let $\pi=\pi_1\otimes\cdots\otimes\pi_k$
be the corresponding discrete series representation of $M$. For
each $i\in\{1,\ldots,k\}$, we fix a $\GL(l_im)$-cover
$(J^{\GL(l_im)},\lambda^{\GL(l_im)})$ of $(J^{\times
l_i},\lambda^{\otimes l_i})$ (as in the proof of
Theorem~\ref{formaldegrees}). Then
\begin{equation} \label{JM}
(J_M,\lambda_M)=(J^{\GL(l_1m)}\times\cdots\times J^{\GL(l_km)},
\lambda^{\GL(l_1m)}\otimes\cdots\otimes\lambda^{\GL(l_km)})
\end{equation}
is a $M$-cover of $(J^{\times e},\lambda^{\otimes e})$.
Then let $(J^G,\lambda^G)$ denote a $G$-cover of
$(J_M,\lambda_M)$ (the existence of which
is guaranteed by \cite[Main Theorem (second version)]{BKss}).

At the same time the partition $(l_1,\ldots,l_k)$ determines the standard Levi
subgroup
\begin{equation} \label{LeviMzero}
M_0=\GL(l_1)\times\cdots\times\GL(l_k)\subset\GL(e,K)=G_0.
\end{equation}

Let $P$ (resp. $P_0$) be the upper-triangular parabolic subgroup of $G$ (resp.
$G_0$) with Levi component $M$ (resp. $M_0$), and unipotent radical denoted by
$N$ (resp. $N_0$).
Let $I$ denote the standard Iwahori subgroup of $G_0$.

\begin{thm} \label{extensiontods}
We have
\[\frac{\vol(J^G\cap N)\cdot \vol(J^G\cap\overline N)}
{\vol_0(I\cap N_0)\cdot \vol_0(I\cap\overline N_0)}=
q^{- \ell(\gamma)f(\sigma^\vee\times\sigma)}
=
\frac{j(\sigma^{\otimes e})}{j_0(1)}
.\]
\end{thm}
\begin{proof}
The second equality follows from our Corollary~\ref{ShahidiIII}.

We will prove the first equality.
Let $U$ denote the unipotent radical of the upper-triangular parabolic subgroup
of $G$ with Levi component $\GL(m)^{\times e}$, and, for $i=1,\ldots,k$,
let $U_i$ denote the unipotent radical of the upper-triangular parabolic
subgroup of $\GL(l_im)$ with Levi component $\GL(m)^{\times l_i}$.
We observe that
\[U=N\times (U\cap M)=N\,\times\,\prod_{i=1}^k U_i.\]
Similarly, let $U_0$ be the unipotent radical of the standard Borel subgroup
of $G_0$, and, for $i=1,\ldots,k$, let $U_{0,i}$ be the unipotent radical of the
standard Borel subgroup of $\GL(l_i,K)$.
We have
\[U_0=N_0\times (U_0\cap M_0)=N_0\,\times\,\prod_{i=1}^k U_{0,i}.\]
It follows from \cite[Proposition~8.5~(i)]{BKtyp} that $(J^G,\lambda^G)$ is also
a $G$-cover of $(J^{\times e},\lambda^{\otimes e})$.
Applying Theorem~\ref{extensiontoeblocks} to $(J^G,U)$ and to
$(J^{\GL(l_im)},U_i)$ for each $i\in\{1,\ldots,k\}$, we obtain
\[\frac{\vol(J^G\cap U)\cdot \vol(J^G\cap\overline U)}
{\vol_0(I\cap U_0)\cdot \vol_0(I\cap\overline U_0)}= q^{-
\frac{e(e-1)}{2}f(\sigma^\vee\times\sigma)} \]
\[\frac{\vol(J^{\GL(l_im)}\cap U_i)\cdot \vol(J^{\GL(l_im)}\cap\overline U_i)}
{\vol_0(I\cap U_{0,i})\cdot \vol_0(I\cap\overline U_{0,i})}=
q^{-\frac{l_i(l_i-1)}{2}f(\sigma^\vee\times\sigma)}.\]
Since $J^G\cap M=J_M$ (by definition of covers), it follows from (\ref{JM}) that
$J^G\cap \GL(l_im)=
J^{\GL(l_im)}$. Then using the fact that
\[\vol(J^G\cap N)=\vol(J^G\cap U)\,\times\,\prod_{i=1}^k\vol(J^{\GL(l_im)}\cap
U_i),\] and the analogous equalities for the others terms, we
obtain
\begin{eqnarray}
\frac{\vol(J^G\cap N)\cdot \vol(J^G\cap\overline N)}
{\vol_0(I\cap N_0)\cdot \vol_0(I\cap\overline N_0)}&=
q^{(-\frac{e(e-1)}{2}+\sum_{i=1}^k\frac{l_i(l_i-1)}{2})f(\sigma^\vee\times\sigma)}\\&
=q^{-\ell(\gamma)f(\sigma^\vee\times\sigma)}.
\end{eqnarray}
\end{proof}

\subsection{Transfer-of-measure}

The following result reduces the case of an arbitrary component $\Omega$
to the one (studied in Corollary~\ref{Iwahori}) of a component (of a possibly
different group $G_0$) which contains the cuspidal pair $(T,1)$.
We give a direct proof which is based on our previous calculations. It is
worth noting that it is also a direct application of \cite[Theorem~4.1]{BHK}.

Let $\Omega=\sigma^e$ be a Bernstein component in $\Omega(\GL(n))$
with single exponent $e$. Let $T$ be the diagonal subgroup of
$G_0=\GL(e,K)$, and let $\Omega_0$ be the Bernstein component in
$\Omega(\GL(e,K))$ which contains the cuspidal pair $(T,1)$. The
components $\Omega$, $\Omega_0$ each have the single exponent $e$,
and we have a homeomorphism of compact Hausdorff spaces
\begin{equation} \label{map}
\Irrt\GL(n,F)_{\Omega} \cong \Irrt \GL(e,K)_{\Omega_0}.
\end{equation}

This homeomorphism is determined by the map
\[
\bigotimes_{i = 1}^k \; \zeta_i^{\val_F \circ \det_F} \otimes \pi_i
\mapsto \bigotimes_{i=1}^k \; (\zeta_i^r)^{\val_K \circ \det_K}
\otimes \St(l_i).\] This formula precisely allows for the fact that
$\pi_i$ has torsion number $r$ and that $\St(l_i)$ has torsion number
$1$. Note that when $\zeta$ is replaced by $\omega \zeta$, where
$\omega$ is an $r$th root of unity, each term remains unaltered.

The equation $r = f(K|F)$ and the standard formula
\[\val_K(y) = f(K|F)^{-1}\, \val_F(N_{K|F}(y))\]
lead to the more invariant formula:
\[ \bigotimes_{i = 1}^k \;
(\chi_i \circ \det\!{}_F) \otimes \pi_i \mapsto \bigotimes_{i=1}^k \;
(\chi_i \circ N_{K|F} \circ \det\!{}_K) \otimes \St(l_i)\] where
$\chi_i$ is an unramified character of $F^{\times}$.

Let $(J^G,\lambda^G)$ be defined as in the previous subsection. It
is a type in $G$ attached to $\Omega$. Recall that $I$ denotes the
standard Iwahori subgroup of $G_0$.

\begin{thm}
Let $d\nu$, $d\nu_0$ respectively denote Plancherel measure on
$\Irrt\GL(n,F)_{\Omega}$, $\Irrt\GL(e,K)_{\Omega_0}$.
We have
\[
\frac{\vol(J^G)}{\dim(\lambda^G)}\,\cdot\,d\nu(\omega) = \vol_0(I)
 \cdot d\nu_0(\omega_0),
\]
where \[\omega = \chi_1\pi_1 \otimes \cdots \otimes \chi_k \pi_k\]
and
\[\omega_0 = (\chi_1 \circ N_{K|F}) \St(l_1) \otimes
\cdots \otimes (\chi_k \circ N_{K|F})\St(l_k).\]
\end{thm}
\begin{proof} We first have to elucidate the canonical measures
$d\omega, d\omega_0$. First, let $M = \GL(n)$, and let $\omega$
have torsion number $r$. Then the map $\Im \, X(M) \to
\mathcal{O}$ is the $r$-fold covering map: $\mathbb{T} \to
\mathbb{T}, z \mapsto z^r$. The map $\Im \, X(M) \to \Im \,
X(A_M)$ sends the map $T \mapsto z^{\val (\det(T))}$ to the map $x
\mapsto z^{val (\det(xI_n))} = (z^{n})^{\val(\det(x))}$ and so
induces the $n$-fold covering map $\mathbb{T} \to \mathbb{T}$. The
canonical measure $\d\omega$ on the orbit $\mathcal{O}$ is the
Haar measure of total mass $n/r$. If $M = \GL(l_1) \times \cdots
\times\GL(l_k)$ and $\omega_j$ has torsion number $r_j$ then the
canonical measure $\d\omega$ on the orbit $\mathcal{O}$ of
$\omega_1 \otimes \cdots \otimes\omega_k$ is the Haar measure of total
mass $l_1 \cdots l_k/ r_1 \cdots r_k$. For the canonical measures
$d\omega, d\omega_0$ we therefore have\[\d\omega = (ml_1 \cdots
ml_k/r^k)\cdot \d\tau = l_1 \cdots l_k \cdot (m^k/r^k) \cdot
\d\tau\] \[\d\omega_0 = l_1 \cdots l_k \cdot \d\tau\] where
$\d\tau$ is the Haar measure on $\mathbb{T}^k$ of total mass $1$.
So, we have
\begin{equation} \label{rapmesures}
\d\omega = (m^k/r^k) \cdot \d\omega_0.
\end{equation}

By Theorem~\ref{PlancherelHCmeasure},
\[
d\nu(\omega)=
q^{\ldee f(\sigma^\vee\times\sigma)}
\,\cdot\,\gamma(G|M)\,\cdot\,d(\omega)\,\cdot\,
\prod
\left|\frac{1 - z_jz_i^{-1}q^{gr}}{1-z_jz_i^{-1}q^{-(g+1)r}}\right|^2
\,\cdot\,
d\omega
\]
and
\[
d\nu_0(\omega_0) = \gamma(G_0|M_0)\,\cdot\,d(\omega_0)\,\cdot\,
\prod \left|\frac{1 -
z_jz_i^{-1}q^{gr}}{1-z_jz_i^{-1}q^{-(g+1)r}}\right|^2\,\cdot\,
d\omega_0.
\]
Hence
\begin{equation} \label{quotientnu}
\frac{d\nu(\omega)}{d\nu_0(\omega_0)}=q^{\ldee f(\sigma^\vee\times\sigma)}
\,\cdot\,\frac{\gamma(G|M)}{\gamma(G_0|M_0)}\,\cdot\,
\frac{d(\omega)}{d(\omega_0)}
\,\cdot\,\frac{d\omega}{d\omega_0}.
\end{equation}

\smallskip

We keep the notation of section~\ref{dc}. It follows from (\ref{volM}),
(\ref{volJ}) that
\begin{equation} \label{volumes}
\vol(J_M)=\vol(J^{\GL(l_1m)})\times\cdots\times \vol(J^{\GL(l_km)}),
\end{equation}
since $J_M=J^{\GL(l_1m)}\times\cdots\times J^{\GL(l_km)}$.
In the same way, we have
\begin{equation} \label{volumesI}
\vol_0(I\cap M_0)=\vol_0(I\cap\GL(l_1m))\times\cdots\times \vol_0(I\cap\GL(l_km)),
\end{equation}

On the other hand, the formula \cite[(7.7.11)]{BK} gives
\[\vol(J^{\GL(l_im)})\,\cdot\,d(\pi_i)=\vol_0(I\cap\GL(l_i,K))\,\cdot\,
\frac{\dim(\lambda^{\GL(l_im)})}{e(E|F)}\,\cdot\,d(\St(l_i))\]

Then (\ref{volumes}), (\ref{volumesI}), (\ref{JM}), and (\ref{productfd}) imply
\begin{equation} \label{extBK}
\vol(J_M)\,\cdot\,d(\omega)= \vol_0(I\cap
M_0)\,\cdot\,\frac{\dim(\lambda_M)}{e(E|F)^k}\,
\cdot\,d(\omega_0).
\end{equation}

Applying (\ref{gammaWalds}) to both $\gamma(G|M)$ and $\gamma(G_0|M_0)$, we
obtain
\begin{equation} \label{quotientgammas}
\frac{\gamma(G|M)}{\gamma(G_0|M_0)}=\frac{\vol(J^G\cap N)
\vol(J^G\cap \overline N)} {\vol_0(I\cap N_0) \vol_0(I\cap
\overline N_0)}\,\cdot\, \frac{\vol_0(I)}{\vol(J^G)}\,\cdot\,
\frac{\vol(J_M)}{\vol_0(I\cap M_0)}.\end{equation}
It then follows from~(\ref{quotientnu}), (\ref{extBK}) and
(\ref{quotientgammas}) that
\[
\frac{d\nu(\omega)}{d\nu_0(\omega_0)}= q^{\ldee
f(\sigma^\vee\times\sigma)} \,\cdot\, \frac{\vol(J^G\cap N)
\vol(J^G\cap \overline N)} {\vol_0(I\cap N_0) \vol_0(I\cap
\overline N_0)}\,\cdot\, \frac{\vol_0(I)}{\vol(J^G)}\,\cdot\,
\frac{\dim(\lambda_M)}{e(E|F)^k}\frac{\d\omega}{\d\omega_0}.\]
Noting that $\dim(\lambda^G)= \dim(\lambda_M)$, and using
equation~(\ref{tors}) and Theorem~\ref{extensiontods}, we have
\[
\frac{\d\nu(\omega)}{\d\nu_0(\omega_0)} =
\frac{\mu_0(I)}{\mu(J^G)} \cdot \dim(\lambda^G) \cdot
\frac{r^k}{m^k} \cdot \frac{\d\omega}{\d\omega_0} =
\frac{\mu_0(I)}{\mu(J^G)} \cdot \dim(\lambda^G), \]
using~(\ref{rapmesures}).

\end{proof}

\section{The central simple algebras case}
Let $D$ be a central division algebra of index $d$ over $F$ and
ring of integers $\integers_D$, and let $A=A(n')$ denote the
algebra of $n'\times n'$ matrices with coefficients in $D$. Then
$A$ is a central simple algebra with centre $F$ of reduced degree
$n=dn'$ and the group of units of $A$ is the group $G'=\GL(n',D)$.
In Theorem 7.2 we will prove a transfer of Plancherel measure
formula for $G'$: this will be deduced from properties of the
Jacquet-Langlands correspondence.
In order to do this, we
will adapt the proof of \cite[(2.5)~p.~88]{AC} to the case when
$F$ is of positive characteristic by using results of
A.~Badulescu.

We use the {\sl standard} normalization of Haar measures, in particular
$\mu_{G'}$ is normalized so that the volume of $\cK'=\GL(n',\integers_D)$ is
$1$.

\subsection{A transfer-of-measure formula}
The aim of this subsection is to prove the transfer-of-measure formula
stated in Theorem~\ref{DivisionAlgebracase}.

An element $x'$ in $G'$ will be called {\sl semisimple} (resp.
{\sl regular semisimple}) if its orbit
$O_{G'}(x')=\left\{yx'y^{-1}\;:\;y\in G'\right\}$ is a closed
subset of $G'$ (resp. if its characteristic polynomial admits only
simple roots in an algebraic closure of $F$). Let $G'_{\rs}$
denote the set of regular semisimple elements in $G'$.

Let $G_{x'}'$ denote the centralizer in $G'$ of $x'$. Then
the group $G_{x'}'$ is unimodular, and the choice of Haar
measures on $G'$ and $G_{x'}'$ induces an invariant measure $dx$ on $G'/G_{x'}'$.
The orbital integral of $f'\in C_c(G')$ at $x'$
is defined as
\begin{equation} \label{intobdef}
\Phi(f',x')=\int_{G'/G_{x'}'}f'(y^{-1}x'y)dy.
\end{equation}
Since the orbit $O_{G'}(x')$ is
closed in $G'$, the integral is absolutely convergent. Indeed, it is a
finite sum, since the restriction of $f'$ to $O_{G'}(x')$ is locally
constant with compact support.
Note that, if $x'\in G'_{\rs}$, then $G_{x'}'$ is a maximal torus in $G'$.

Orbital integrals have a local expansion, due to Shalika \cite{Sha}, which we will
now recall. If $O'$ is a unipotent orbit in $G'$, let $\Lambda_{O'}$ denote
the distribution given by integration over the orbit $O'$.
There exist functions $\Gamma_{O'}^{G'}\colon G'_{\rs}\to\RR$ (the {\sl Shalika
germs}) indexed by unipotent orbits of $G'$ with the following property:
\begin{equation} \label{germexpansion}
\Phi(f',x')=\sum_{O'}\,\Gamma_{O'}^{G'}(x')\,\cdot\,\Lambda_{O'}(f'),\end{equation}
for $x'\in G'_{\rs}$ sufficiently close to the identity.
Observe that $\Lambda_1=f'(1)$.

Harish-Chandra proved that the germ $\Gamma_1^{G'}$ associated to the trivial
unipotent orbit is constant, and Rogawski \cite{R} has determined its value
assuming the characteristic of $F$ to be zero:
\begin{equation} \label{Rogawski}
\Gamma_1^{G'}=\frac{(-1)^{n-n'}}{d(\St_{G'})}.\end{equation}
The equality (\ref{Rogawski}) is still valid in the case when $F$ is of
positive characteristic.
Indeed, let $F$ be of positive characteristic and let $E$ be a field of zero
characteristic sufficiently close to $F$, that is, such that there exists a ring
isomorphism from $\integers_F/\varpi^l\integers_F$ to
$\integers_E/\varpi^l\integers_E$, for some sufficiently big
integer $l\ge 1$. Let $D_E$ be a central division algebra over $E$
with the same index $d$. Then by \cite[Lemma~3.8]{Badintorb} the
lifts $f'_E$ of $f'$ to $G'_E=\GL(m,D_E)$ (resp. $f_E$ of $f$ to
$G_E=\GL(n,E)$) also satisfy $f_E\leftrightarrow (-1)^{n-n'}
f'_E$. On the other hand, $f'_E(1)=f'(1)$, independently of $m$:
since the way to lift $f'$ to $f'_E$ consists in cutting the group
$G'$ into compact open subsets on which $f'$ is constant, in
associating to these subsets compact open subsets in $G'_E$, and
assigning to \emph{these} subsets the same constants in order to
define $f'_E$; but the compact open subset of $G'$ containing $1$
corresponds to the compact open subset in $G'_E$ containing $1$.

\smallskip

\smallskip

If $\pi$ is a smooth representation of $G$ or $G'$ with finite
length, we will denote by $\theta_\pi$ its character.

\begin{thm} \label{JLCor}
{\rm The Jacquet-Langlands correspondence \cite{DKV}, \cite{BadJL}.}
There exists a bijection
\[\JL\colon \cE_2(G')\to \cE_2(G)\]
such that for each $\pi'\in\cE_2(G')$:
\begin{equation} \label{JLC}
\theta_{\pi'}(x')=(-1)^{n-n'}\,\theta_{\JL(\pi')}(x),
\end{equation}
for any $(x,x')\in G\times G'$ such that $x\leftrightarrow x'$.
\end{thm}

Recall that $A=A(n')$ denotes the algebra of $n'\times n'$
matrices with coefficients in $D$. Let $\Nrd_{A|F}\colon A\to F$
denote the reduced norm of $A$ over $F$ as defined in
\cite[\S~12.3, p.~142]{B2}. We shall view the reduced norm
$\Nrd_{A|F}$ as a homomorphism from $G'$ to $F^\times$.

If $\eta$ is a quasicharacter of $F^{\times}$ then we will write
\[
\eta \pi' = (\eta\circ\Nrd_{A|F})\otimes\pi'.
\]
If $\eta$ is an unramified quasicharacter then we will refer to
$\eta \pi'$ as an \emph{unramified twist} of $\pi'$.

Each representation $\pi'$ of $G'$ has a {\it torsion number}: the
order of the cyclic group of all those unramified characters
$\eta$ of $F^\times$ for which
\[\eta \pi' \cong \pi'.\]
The Jacquet-Langlands correspondence has the property that
\begin{equation} \label{torsioncJL}
\eta (\JL(\pi'))= \JL(\eta \pi'),
\end{equation}
for any square integrable representation $\pi'$ of $G'$ and any
(unitary) character $\eta$ of $F^\times$ (see
\cite[(4)~p.~35]{DKV}). It follows that the torsion number of
$\pi'$ is equal to that of $\JL(\pi')$.

For each Levi subgroup
$M=\GL(n_1,F)\times\cdots\times\GL(n_k,F)$ of $G$ such that $d$ does
not divide $n_i$ for some $i\in\{1,\ldots,k\}$, we have
\[\theta^{G}_{\omega}(f)=0, \;\;\text{ for any
$\omega\in\cE_2(M)$}\]
(see the beginning of \cite[\S 3]{Badintorb} and the proof of
\cite[Lem.~3.3]{Badintorb}).

We consider now a Levi subgroup $M$ of the form
$M=\GL(dn'_1,F)\times\cdots\times\GL(dn'_k,F)$, and define
$M'=\GL(n_1',D)\times\cdots\times\GL(n_k',D)$ (a Levi subgroup of
$G'$): $M$ is the \emph{transfer} of $M'$. The Jacquet-Langlands
correspondence induces a bijection
$\JL\colon\cE_2(M')\to\cE_2(M)$, by setting
\[\JL(\omega_1'\otimes\cdots\otimes\omega_k')
=\JL(\omega'_1)\otimes\cdots\JL(\omega_k').\] For any
$\omega\in\cE_2(M)$, there exists $\omega'\in\cE_2(M')$ such that
$\omega=\JL(\omega')$.

Let $\Omega^\temp(G')$, $\Omega^\temp(G)$ denote the Harish-Chandra
parameter space of $G', G$.  Each point in $\Omega^\temp(G')$ is a
$G'$-conjugacy class of discrete-series pairs $(M',\omega')$ with
$\omega' \in E_2(M')$.  The topology on $\Omega^\temp(G')$ is
determined by the unramified unitary twists: then $\Omega^\temp(G')$
is a locally compact Hausdorff space.  The map
\[
(M',\omega') \mapsto (M,\JL(\omega')),
\]
where $M$ is the transfer of $M'$, secures an \emph{injective} map
\[
\JL\colon\Omega^\temp(G') \to \Omega^\temp(G). \]  We will write
$Y = \JL(\Omega^t(G'))$.  Since the $\JL$-map respects unramified
unitary twists, we obtain a homeomorphism of $\Omega^\temp(G')$
onto its image:
\[
\JL\colon \Omega^\temp(G') \cong Y \subset \Omega^\temp(G).
\]



\begin{thm} \label{DivisionAlgebracase} {\rm Transfer of Plancherel
measure.} Let $G' = \GL(n',D), G = \GL(n,F)$ with $n = dn'$.  Let
$\nu', \nu$ denote the Plancherel measure for $G',G$, each with
the standard normalization of Haar measure on $G',G$. Then we have
\[
d\nu'(\omega') = \lambda(D/F) \cdot d\nu(\JL(\omega'))
\]
where
\[
\lambda(D/F) = \prod(q^m - 1)^{-1}\] the product taken over all
$m$ such that $1 \leq m \leq n-1, m \neq 0 \, \mod \, d$.

\end{thm}
\begin{proof}
If $x\in G$ and $x'\in G'$, we will write
$x\leftrightarrow x'$ if $x$, $x'$ are regular semisimple and have the same
characteristic polynomial. If $x\in G$, we will say that {\sl $x$ can be
transferred} if there exists $x'\in G'$ such that $x\leftrightarrow x'$.

Let $f'\in C_c(G')$. Then, by \cite[Th.~3.2.]{Badintorb}, there exists
$f\in C_c(G)$ such that
\[
\Phi(f,x)=
\begin{cases}
(-1)^{n-n'}\cdot\Phi(f',x')&\text{for each $x'\in G'$ such that $x\leftrightarrow x'$,}\\
0&\text{ if $x$ cannot be transferred,}
\end{cases}
\]
for any $x\in G_{\rs}$.

It then follows from the germ expansion~(\ref{germexpansion}) that
\[f'(1)\,\cdot\,\Gamma_1^{G'}=
(-1)^{n-n'}\,\cdot\,f(1)\,\cdot\,\Gamma_1^{G},\]
that is, using~(\ref{Rogawski}),
\begin{equation} \label{f1}
\frac{f'(1)}{d(\St_{G'})}=\frac{f(1)}{d(\St_{G})}.
\end{equation} We recall that
$\theta^G_{\omega}(f) = 0$ on the complement of $Y$ in
$\Omega^t(G)$. Next, we use equation (\ref{f1}), and apply twice
the
 Harish-Chandra Plancherel theorem, first for $G'$, then for $G$. We obtain
\begin{eqnarray}\label{formule10}
\int\theta^{G'}_{\omega'}(f')\,d\nu'(\omega')& = &
f'(1)\nonumber\\& = & d(\St_{G'})\cdot d(\St_{G})^{-1} \cdot f(1)
\nonumber\\&=& d(\St_{G'})\cdot d(\St_{G})^{-1} \cdot
\int\theta^{G}_{\omega}(f)\,d\nu(\omega)\nonumber\\
&=&d(\St_{G'})\cdot d(\St_{G})^{-1} \cdot
\int\theta^{G}_{\omega}(f)\,d\nu|_{Y}(\omega), \label{HCequality}
\end{eqnarray}
for all $f'\in C_c(G')$.

We recall that the parameter space $\Omega^\temp(G')$ is the
\emph{domain} of the Plancherel measure $\nu'$.

By the refinement of the trace Paley-Wiener theorem due to
Badulescu \cite[Lemma~3.4]{Badintorb} we have
\[
\{\omega' \mapsto \theta^{G'}_{\omega'}(f^{\prime\vee}) : f' \in C_c(G'),
\omega' \in \Omega^\temp(G')\} = L(\Omega^\temp(G')),
\]
where $L(\Omega^\temp(G'))$ is the space of compactly supported
functions on $\Omega^\temp(G')$ which, upon restriction to each
connected component (a quotient of a compact torus $\mathbb{T'}^k$
by a product of symmetric groups), are Laurent polynomials in the
co-ordinates $(z_1, z_2, \ldots,z_k)$.

Now $L(\Omega^\temp(G'))$ is a dense subspace of
$C_0(\Omega^\temp(G'))$, the continuous complex-valued functions
on $\Omega^\temp(G')$ which vanish at infinity. On the other hand,
it follows from \cite[Prop.~3.6]{Badintorb} that
\begin{equation} \label{traces}
\theta^{G'}_{\omega'}(f')=
\theta^{G}_{\JL(\omega')}(f),\;\;\text{for any
$\omega'\in\cE_2(M')$.}
\end{equation}

Equation (\ref{HCequality}) therefore provides us with two Radon
measures (continuous linear functionals) which agree on a
\emph{dense subspace} of $C_0(\Omega^\temp(G'))$. Therefore the
measures are equal:

\begin{equation}
\label{Formule3} d\nu'(\omega') = d(\St_{G'})\cdot d(\St_{G})^{-1}
\cdot d\nu|_{Y}(\omega)
\end{equation}

At this point, we have to elucidate a normalization issue. Let $K'
= \GL(n',\mathfrak{o}_D)$. The group $A_{G'}$ by definition is the
$F$-split component of the centre of $G'$ and can be identified
with $F^\times$. As in section (6.2), we have
$F^{\times}K'/F^{\times} = K'/K' \cap F^{\times} =
K'/\mathfrak{o}_F^{\times}$. But the Haar measure on $A_{G'}$ has,
as in \cite[p.240]{Wald}, the standard normalization $\mes(K' \cap
A_{G'}) = 1$, \ie $\mes(\mathfrak{o}^{\times}_F) = 1$. Since
$\mes(K') = 1$, we have $\mes(F^{\times}K'/F^{\times}) = 1$. It
follows (see for instance \cite[3.7]{SZ}) that the formal degree
of the Steinberg representation $\St_{G'}$ is given by
\[d(\St_{G'})= \frac{1}{n} \prod_{j=1}^{n'-1}(q^{dj}-1)\]

We then have
\begin{equation}
\label{Formule4} d\nu'(\omega') = \lambda(D/F) \cdot d\nu(\omega)
\end{equation}
where
\[
\lambda(D/F) = (q^d-1)(q^{2d}-1) \cdots
(q^{(n'-1)d}-1)(q-1)^{-1}(q^2-1)^{-1} \cdots (q^{n-1}-1)^{-1},\]
so that
\begin{equation}
\label{Formule5} \lambda(D/F) = \prod(q^m-1)^{-1}
\end{equation}
the product taken over all $m$ such that $1 \leq m \leq n-1$, $m
\neq 0 \mod \,d$.
\end{proof}
This result may be expressed as follows
\begin{thm} Let
$(\Omega^\temp{G'}, \mathcal{B'}, \nu')$ be the measure space
determined by the Plancherel measure $\nu'$, let $(Y,
\mathcal{B},\lambda(D/F) \cdot \nu|_{Y})$ be the measure space
determined by the restriction of $\lambda(D/F) \cdot \nu$ to $Y =
\JL(\Omega^\temp(G') \subset \Omega^\temp(G)$. Then these two measure
spaces are isomorphic:
\[(\Omega^\temp{G'}, \mathcal{B'}, \nu') \cong
(Y, \mathcal{B},\lambda(D/F) \cdot \nu|_{Y})\]
\end{thm}

Anne-Marie Aubert, Institut de Math\'ematiques de Jussieu,
U.M.R. 7586, Universit\'e
Pierre et Marie Curie, F-75252 Paris Cedex 05,
France\\
Email: aubert@math.jussieu.fr

Roger Plymen, Mathematics Department, Manchester University, M13 9PL, U.K.\\
Email: plymen@manchester.ac.uk
\end{document}